\documentclass[11pt]{article}
\usepackage{amssymb,amsmath}
\usepackage{color}
\usepackage{graphics}
\usepackage{graphicx}  
\numberwithin{equation}{section}
\newtheorem{theorem}{Theorem}[section]                   %theory and  corollary
\newtheorem{proposition}{Proposition}[section]
\newtheorem{lemma}{Lemma}[section]

\newtheorem{remark}{Remark}[section]

\newcommand{\be}{\begin{equation}}
\newcommand{\ee}{\end{equation}}

\newcommand{\into}{{\int_{\Omega}}}

\newcommand{\la}{\lambda}

\def\la{\lambda}

\begin{document}
	\title{ Slow decay and turnpike for infinite-horizon
		hyperbolic LQ problems\thanks{The first author  was supported by the Natural Science Foundation of 
			China grant NSFC-62073236. The second author was supported by the European Union's Horizon 2020 research and innovation programme under the Marie Sklodowska-Curie grant agreement No.765579-ConFlex, from the Alexander von Humboldt-Professorship program, the European Research Council (ERC) under the European Union's Horizon 2020 research and innovation programme (grant agreement NO. 694126-DyCon), the Transregio 154 Project ``Mathematical Modeling, Simulation and Optimization Using the Example of Gas Networks” of the German DFG, and grant MTM2017-92996-C2- 1-R COSNET of MINECO (Spain). }}
	
	\author{Zhong-Jie Han \footnote{{School of Mathematics,  Tianjin University, Tianjin 300354, China.  (Email address:  zjhan@tju.edu.cn).
		}}
		\quad Enrique Zuazua \footnote{
			Chair for Dynamics, Control and Numerics-Alexander von Humboldt-Professorship,	Department of Data Science,
			Friedrich-Alexander-Universitat Erlangen-Nurnberg,
			91058 Erlangen, Germany. (Email address: enrique.zuazua@fau.de).}
		\footnote{Chair of Computational Mathematics, 
			Fundacion Deusto, Av. de las Universidades, 24,
			48007 Bilbao, Basque Country, Spain.}
		\footnote{Departamento de Matem\'{a}ticas,
			Universidad Aut\'{o}noma de Madrid, 28049 Madrid, Spain.}}
	\maketitle
	
	\begin{abstract}
		This paper is devoted to analysing the explicit slow decay rate and turnpike in the infinite-horizon linear quadratic optimal control problems for hyperbolic systems. Assume that some weak observability or controllability are satisfied, by which, the lower and upper bounds of the corresponding algebraic Riccati operator are estimated, respectively. Then based on these two bounds, the explicit slow decay rate of the closed-loop system with Riccati-based optimal feedback control  is obtained. The averaged turnpike property for this problem is also further discussed. We then apply these results to the  LQ optimal control problems constraint to networks of one-dimensional wave equations and also some multi-dimensional ones with local controls which lack of  GCC(Geometric Control Condition). 
	\end{abstract}
	
	\noindent {\bf Key words:}\hspace{2mm}  \vspace{3mm} 
		Optimal control problems,  Riccati operator, slow decay rate, weak controllability and observability, turnpike property.

	\noindent {\bf AMS subject classifications.}\hspace{2mm}  \vspace{3mm} 
		49J20, 49K20, 93C20,  49N05.

%	\pagestyle{myheadings}
%	\thispagestyle{plain}
%	\markboth{Z.J. HAN AND E. ZUAZUA}{Slow decay and turnpike for infinite-horizon  hyperbolic LQ problems}
%\title{ Slow decay rate for infinite horizon optimal problem}
%% \thanks{This research is supported by the Natural Science
%%	Foundation of China grant NSFC-61573252.} 
%
%%\date{Jun. 1, 2004.}
%\author{Zhong-Jie Han$^1$, \qquad Enrique Zuazua$^2$} 
%School of Mathematics,  Tianjin University, Tianjin 300354, China. 
%
%%zjhan@tju.edu.cn\qquad gqxu@tju.edu.cn
%%Department of Auto
%\maketitle
%\begin{abstract}
%This work studies the slow turnpike property for the networks of wave equations.
%
%First, we prove the convergence of the average with a rate under smooth initial condition.
%
%Second, we further estimate the rate for slow turnpike property for this network.
%
%
%\end{abstract}

%\noindent {\bf Key words:}\hspace{2mm}  \vspace{3mm} 

\section{Introduction}
The object of this paper  is devoted to discussing the large time behaviour and turnpike property in  infinite-time  linear quadratic(LQ) optimal control problems under weak  controllability and observability hypothesises.  Specifically,  we will discuss the relationship between the bounds of the corresponding algebraic Riccati operator and the weak controllability and observability properties, and based on which  we  identify the explicit slow decay rate of the closed-loop system with the Riccati-based optimal feedback control.  Moreover, under  weak controllability  and observability hypothesises,  we further discuss that how the optimal control and trajectories of the LQ optimal control problems  converge to the corresponding stationary optimal control and state, that is the so-called turnpike property.
   
   LQ optimal control problems have been studied extensively in recent thirties years, see
 \cite{willems} for finite dimensional systems,  \cite{BDDM},  \cite{datko}, \cite{morris}, \cite{Li} and \cite{sontag} for the infinite dimensional systems with bounded input or output operator, and \cite{DD}, \cite{sammon}, \cite{Russell} for the ones with unbounded input or output operator.
  
It is known that the solution of LQ
optimal control problems can be constructed as a feedback form based on solving its corresponding Riccati equations(see  \cite{curtain}, \cite{sontag}). In other words,
the solution to (algebraic) Riccati equation, which is called (algebraic)  Riccati operator, is corresponding to the  finite-time (or infinite-time) LQ optimal feedback control. Especially, the algebraic Riccati operator is usually used to stabilize the system. For instance, Aksikas et. al. in \cite{LQapp} designed the LQ feedback control to stabilize a class of hyperbolic PDE systems exponentially, by solving the matrix Riccati differential equation. 
Porretta and Zuazua  in \cite{zuazua2} and \cite{zuazua3} proved the exponential decay of the closed-loop system with Riccati-based optimal feedback control under the assumptions of exact observability of $(A, C)$ and exact controllability of $(A, B)$. Indeed, based on the exact observability  of $B^*$ and $C$, they obtained that $\langle\hat{E}x, x\rangle$ is a strict Lyapunov function for the system, where $\hat{E}$ is the algebraic Riccati operator and $x$ is the state of the system in the  Hilbert state space $H$. Meanwhile, it is also proposed in \cite{zuazua2} that the exponential-type point-wise turnpike property can be further proved based on the exponential decay rate of the closed-loop system with Riccati-based optimal feedback controls (see  \cite{turnpike1}, \cite{turnpike2} and \cite{turnpike3}).   
 
Note that in  the previous results,  both the stabilization of infinite-time LQ optimal control problems and turnpike properties  are all discussed under the assumptions of exact controllability and observability. 
In  \cite{zuazua2},  a kind of slow turpike in average(Logarithmic-type)  was concluded for multi-dimensional wave equation without GCC(Geometric Control Condition), which inspired us to give a complete analysis on  the turnpike properties for the  LQ optimal control problems if lacking of exact controllability or observability. 
It is obvious that  some weaker controllability or observability are still necessary so as to guarantee  not only the existence of the solutions to  infinite-time LQ optimal control problems  but also the feedback stabilization of the  system.  Thus, to do this,  some  hypothesises on weak-type observability estimates  are chosen as given in the next section.   

In this work, we shall first  address the following problem:   

Q1. To which extent  do  weak observability or controllability determine the decay rates of the energy for the systems with Riccati-based optimal feedback controls?  

We find that if the existence of the solution to the LQ optimal control problems is guaranteed, then the lower bound of $\langle\hat{E}x, x\rangle$ is totally determined by the observability property of output operator $C$, while the upper bound is totally dependent on the controllability property of the input operator $B$. Specifically, on one hand, when $C$  is weakly observable, the lower bound can be obtained respect to an weak norm of $x$, which is estimated completely based on the extent of the weak observability. On the other hand, when  $B$  is weakly controllable, then the upper bound of $\langle\hat{E}x, x\rangle$ can be estimated respect to a strong norm of $x$,  totally by the extent of the weak controllability. Based on these properties of the algebraic Riccati operator, we deduce the explicit slow decay rate of the closed-loop system related to the infinite-time LQ optimal control problems. 

The solving of Q1 is one key step  to  further discuss the turnpike properties of the LQ optimal control problems(see \cite{zuazua2}, \cite{zuazua3}), that is the following problem  under consideration:

Q2. To which extent do weak observability or controllability lead to turnpike properties of the LQ optimal control problems?

Based on the weak observability of $(A,B^*)$ and $(A,C)$,  we  show that the averaged turnpike property of such  problems holds under certain conditions on the initial state, the stationary optimal state and its dual.
It is worth mentioning that  in \cite{zuazua2}, under the exact controllability and observability, 
 the exponential decay of the closed-loop system with Riccati-based optimal control always holds, based on which
 the exponential-type point-wise turnpike property for the LQ problems can be proved.  Thus, under weak controllability and observability, note that  the slow decay rate of the closed-loop system can be estimated in our this work, and then the slow(polynomial-type etc.) point-wise turnpike property
should be also reasonable to be expected.  

However,  the (slow) point-wise turnpike property of the LQ optimal control problems under weak controllability or observability hypothesises  is still  an open problem and  worth further investigating in  future. In fact,   the non-uniform slow decay rates of the closed-loop system inevitably caused by the weak observability of $(A,B^*)$ and $(A,C)$ are always dependent on the regularity of the initial states. This  makes  the energy estimates for (slow) point-wise turnpike property become a tough issue to be tackled.  This is different from the case under exact controllability and observability,  in which the uniform exponential decay of the closed-loop systems always holds for all initial states and  the corresponding energy estimates can be easily carried out.

The rest of this paper is given as follows. In Section \ref{main}, the preliminary and main results of this paper are presented. The bounds estimates for the algebraic Riccati operator and the explicit slow decay rate are given under weak controllability and observability hypothesises. The averaged turnpike property is also further  presented.   In Section \ref{333}, we  prove the main results given in this work.   Section \ref{examples} is devoted to presenting  some  examples on the slow decay rates and turnpike for  some kinds of hyperbolic systems(networks of wave equations and multi-dimensional ones)  without GCC. Finally, in Section \ref{open} a conclusion and future work  are given.

\section{Preliminary and main results} \label{main}
This section is devoted to problem formulation and further presenting the decay result and turnpike property of the LQ optimal control problems under the weak assumptions of the  controllability and observability.
\subsection{Problem formulation}
Similar to the abstract frame setting for control
 systems in \cite{TW},
assume that  $A: \mathcal{D}(A)\to H$ is a self-adjoint, strictly positive operator with compact resolvent. Thus, $A$ is diagonalizable.  Hence, if the eigenvalues  of $A$ are  given as $(\lambda_n^2)_{n\geq 1}$,  its corresponding eigenvectors $(\phi_n)_{n\geq 1}$ forms an  orthonormal basis in $H$. 
Moreover, $$D(A):=\{\varphi \in H| \sum_{n\geq 1} \lambda_n^{4} |\langle \varphi, \phi_k \rangle|^2<\infty \}.$$

Define the space $X_\beta$ as follows.
$$X_\beta:=\{\varphi \in H| \sum_{n\geq 1} \lambda_n^{4\beta} |\langle \varphi, \phi_k \rangle|^2<\infty \}$$
with inner product
$$\langle \varphi, \psi \rangle_\beta =\langle A^{\beta} \varphi, A^{\beta} \psi \rangle$$
where $$A^\beta \varphi:=\sum_{n\geq 1}\limits \lambda_n^{2\beta}\langle z, \phi_n\rangle\phi_n.$$
Let us consider the following control system: 
\begin{equation}\label{001}
\left\{
\begin{array}{l}
w_{tt}+A w=B u(t),\\
w(0)=w_0,\quad w_t(0)=w_1.
\end{array}
\right.
\end{equation}
%As before, we assume that $A$ satisfies \rife{coerc}.
%Assume that $X$, $H$ satisfy $X\subset H\subset X'$.
%  and satisfy
% \begin{equation}
% \exists \la\,,\mu >0:   \quad \langle Ax, x\rangle_{X',X}  + \mu |x|^2_H\geq \la\, \|x\|^2_X\qquad \forall x\in X.
% \end{equation}
In (\ref{001}), assume that $B\in {\mathcal L}(U, H)$  and {$C \in {\mathcal L}(X_{{1}\!/\!{2}},V)$}, where $U$, $V$ are Hilbert spaces. 
We know that system (\ref{001}) is well-posed with input space $U$ in Hilbert state space $\mathcal{H}:=X_{{1}\!/\!{2}}\times H$ with norm 
$$\|(w,v)\|_{\mathcal{H}}^2=\langle(w,v), (w,v)\rangle_{\mathcal{H}}=\langle A^{\frac{1}{2}} w, A^{\frac{1}{2}} w\rangle+ \langle v,v \rangle.$$
Thus,  for $(w_0, w_1)\in \mathcal{H}$ and $u\in L^2([0,\infty), U)$, system (\ref{001}) admits a unique solution $w\in C([0,\infty), X_{{1}\!/\!{2}})\cap C^1([0,\infty); H)$.  {For simplicity, we choose $U=V=H$}.

Set  $\mathcal{A}:=\left[ \begin{array}{cc} 0&I \cr -A& 0\end{array} \right]$ with $\mathcal{D}(\mathcal{A})=X_{1}\times X_{{1}\!/\!{2}}$, and $\mathcal{B}=\left[\begin{array}{c} 0\\ B \end{array}\right]$. Then system (\ref{001}) also can be written as follows:
\begin{equation}\label{001+}
\left\{
\begin{array}{l}
\frac{d W}{d t}=\mathcal{A} W+\mathcal{B}(0, u(t))^T,\\
W(0)=(w_0, w_1)^T,
\end{array}
\right.
\end{equation}
where  $W= (w, w_t)^T$.

 It is well-known that $\mathcal{A}$ generates a $C_0$ contraction semigroup on $\mathcal{H}$.
Define the interpolation space
$\mathcal{D}(\mathcal{A}^s)=X_{{(s+1)}\!/\!{2}}\times X_{{s}\!/\!{2}}$ and so its dual one
$\mathcal{D}{(\mathcal{A}^s)'}=X_{-{s}\!/\!{2}}\times X_{-{(s+1)}\!/\!{2}}$. Thus, we get 
\begin{equation}\label{da1}
 \|(w_0, w_1)\|_{\mathcal{D}(\mathcal{A}^{s})}^2=\sum_{j\geq 1}\limits
\lambda_n^{2s}(\lambda_n^2a_n^2+b_n^2),
\end{equation}
and
\begin{equation}\label{da2}
 \|(w_0, w_1)\|_{\mathcal{D}(\mathcal{A}^{s})'}^2=\sum_{j\geq 1}\limits
\lambda_n^{-2(s+1)}(\lambda_n^2a_n^2+b_n^2),
\end{equation}
where
 $a_n,\; b_n$ are the Fourier's coefficients and
\begin{equation}\label{fourier}
w_0=\sum_{n\geq 1}a_n \phi_n(x),\quad w_1=\sum_{n\geq 1}b_n \phi_n(x).
\end{equation}

Assume that the operator $B$ and $C$ satisfy the following two estimates.

{\bf (H1).}  $(A,B^*)$ is weakly observable, that is, for the system 
\begin{equation}\label{02}
\left\{
\begin{array}{l}
w_{tt}+Aw=0,\\
w(0)=w_0,\quad w_t(0)=w_1,
\end{array}
\right.
\end{equation}
there exist positive constants $T_0$ and $c$ such that 
\begin{equation}\label{03}
c\int_0^{T_0}\|B^*w_t\|_{U'}^2dt\geq\| (w_0, w_1)\|_{X_{{1}\!/\!{2}-{1}\!/\!{(2\varrho)}}\times X_{-{1}\!/\!{(2\varrho)}}}^2=\sum_{j\geq 1}\limits
\lambda_n^{-\frac{2}{\varrho}}(\lambda_n^2a_n^2+b_n^2),
\end{equation}
where $\varrho>0$ is some constant, $\lambda_n,\; a_n,\; b_n,\; n=1,2,...$ are given as in (\ref{fourier}). 

{\bf (H2).} $(A,C)$ is weakly observable, that is, there exist positive constants $T_1$ and $c$ such that 
{\begin{equation}\label{obserC}
c\int_0^{T_1}\|C{w}\|_V^2dt\geq\| (w_0, w_1)\|_{X_{{1}\!/\!{2}-{1}\!/\!{(2\eta)}}\times X_{-{1}\!/\!{(2\eta)}}}^2{=\sum_{j\geq 1}\limits
 \lambda_n^{-\frac{2}{\eta}}(\lambda_n^2a_n^2+b_n^2),}
\end{equation}}
  where $\eta>0$ is some constant, $\lambda_n,\; a_n,\; b_n,\; n=1,2,...$ are  also given as in (\ref{fourier}). 

%  Due to $C\in \mathcal{L}(H, H)$, this implies that 
%{\begin{equation}\label{obserC1}
%	c\int_0^{T_1}\|CA^{\frac{1}{2}}w\|_V^2dt\geq\| (w_0, w_1)\|_{X_{\frac{1}{2}-\frac{1}{2\eta}}\times X_{-\frac{1}{2\eta}}}^2{=\sum_{j\geq 1}\limits
%		\lambda_n^{-\frac{2}{\eta}}(\lambda_n^2a_n^2+b_n^2),}
%	\end{equation}}

%\begin{remark}
%	It should be noted that in the above hypothesises, if $\varrho$ or $\eta$ goes to infinity, it implies the exact  controllability of $(A,B)$
%	or observability of $(A,C)$.
%\end{remark}

Based on (H1), together with the  semigroup theory, we can get the following slow decay rate for the system with the collocated feedback controls (see \cite{ammari-tucsnak}).

\begin{lemma}\label{le1-11} Assume that (H1) is fulfilled. 
Under the feedback control law 
	\begin{equation}\label{002}
	u(t)=-B^* w_t,
	\end{equation}
	 it holds that  for all $t>0$, the solution to the closed-loop system (\ref{001}) decays polynomially for any $(w_0, w_1)\in {\mathcal{D}(\mathcal{A}^k)},\; k>0$, that is 
	\begin{equation}\label{3}
	\|(w, w_t)\|^2_{\mathcal{H}}\leq c_1 (t+1)^{-k \varrho }\|(w_0, w_1)\|_{\mathcal{D}(\mathcal{A}^k)}^2,\; k>0,
	\end{equation}
	where $\varrho>0$ is the same as in { (H1)} and $c_1$ is a constant independent of initial data.
\end{lemma}

Thus,  by multiplying (\ref{001}) with $w_t$, we have  
\begin{eqnarray}\label{20002+}
\int_0^T ||B^*w_t||_{H}^2dt&=&	-\frac{1}{2}\|(w(T), w_t(T))\|^2_{\mathcal{H}}+\frac{1}{2}\|(w_0, w_1)\|^2_{\mathcal{H}}\cr
&\leq&  \frac{1}{2}\|(w_0, w_1)\|^2_{\mathcal{H}}+\frac{1}{2}c_1 (T+1)^{-{k} {\varrho} }\|(w_0, w_1)\|_{\mathcal{D}(\mathcal{A}^{{k}})}^2.
\end{eqnarray}

{ If using  the transformation $\phi=\int_0^t wdt+\phi(0)$ in (\ref{02}), it holds that
	$$
	\left\{
	\begin{array}{l}
	\phi_{tt}+A\phi=0,\\
\phi(0)=A^{-1}w_1,\quad \phi_t(0)=w_0.
	\end{array}
	\right.
$$

By (H2), we get	
\begin{equation}\label{obse2}
c\int_0^{T_1}\|C\phi_t\|_V^2dt\geq\| (\phi(0), \phi_t(0))\|_{X_{1-{1}\!/\!{(2\eta)}}\times X_{{1}\!/\!{2}-{1}\!/\!{(2\eta)}}}^2.
\end{equation}
	
	Then, similar to Lemma \ref{le1-11},  we obtain the following result.

 \begin{lemma}\label{le1-12}
 	Assume that (H2) is fulfilled. Then for all $0<t<T$ and $(\phi_0^T, \phi_1^T)\in {\mathcal{D}(\mathcal{A}^{k})},\; k\geq 1$, the solution to the following backward closed-loop system
 \begin{equation}\label{02+}
 \left\{
 \begin{array}{l}
 \phi_{tt}+A \phi=C^*C\phi_t,\\
 \phi(T)=\phi_0^T\in X_1,\quad \phi_t(T)=\phi_1^T\in X_{{1}\!/\!{2}},
 \end{array}
 \right.
 \end{equation}
 satisfies
 	\begin{equation}\label{3+}
 	\|(\phi, \phi_t)\|^2_{\mathcal{H}}\leq c_1 (T-t+1)^{-k \eta }\|(\phi_0^T, \phi_1^T)\|_{\mathcal{D}(\mathcal{A}^{k})}^2,\; k\geq 1,
 	\end{equation}
 	where $\eta>0$ is the same as in { (H2)} and $c_1$ is a constant independent of terminal data.
 \end{lemma}
}

Thus,  by multiplying (\ref{02+}) with $\phi_t$, we have  
\begin{eqnarray}\label{214}
\int_0^T ||C\phi_t||_{H}^2dt&=&	\frac{1}{2}\|(\phi_0^T, \phi_1^T)\|^2_{\mathcal{H}}-\frac{1}{2}\|(\phi(0), \phi_t(0))\|^2_{\mathcal{H}}\cr
&\leq&  \frac{1}{2}\|(\phi_0^T, \phi_1^T)\|^2_{\mathcal{H}}+\frac{1}{2}c_1 (T+1)^{-{k} \eta }\|(\phi_0^T, \phi_1^T)\|_{\mathcal{D}(\mathcal{A}^{{k}})}^2,\; k\geq 1.
\end{eqnarray}
%{\bf Hypothesis 1a.}  Set the state space $\mathcal{H}=H^1\times L^2$. Under the feedback law 
%\begin{equation}\label{002}
%u(t)=-B^* y_t,
%\end{equation}
%the closed-loop system (\ref{001}) decay polynomially, that is 
%\begin{equation}\label{3}
%\|y, y_t)\|^2_{\mathcal{H}}\leq C (t+1)^{-k \varrho }\|(y_0, y_1)\|_{\mathcal{D}(\mathcal{A}^k)}^2,
%\end{equation}
%where $\varrho>0$ is some constant.

We consider the following quadratic performance index associated with the control system (\ref{001}).

\begin{equation}\label{opt}
\min J_t^T(u)=\frac{1}{2}\int_t^T [\|u(t)\|^2_U
+\|Cw(t)\|_V^2]dt,\; u\in L^2(t,T; U) .
\end{equation}
The corresponding  OS (Optimality System)  is given as follows:
\begin{equation}\label{optimal system1}
\left\{
\begin{array}{l}
w_{tt}^T + Aw^T =   Bu^T, \quad t\leq s\leq T, \\
w(t)=w_0, \quad w_t(t)=w_1,\\
u^T= -B^*{ p^T}, \\
p^T_{tt} + Ap^T=C^*C w^T, \quad t\leq s\leq T, \\
p^T(T)=p_t^T(T)=0.
\end{array}\right.
\end{equation}
	\begin{lemma}
		Assume that (H2) is fulfilled. Then there exists a unique linear operator $ \mathcal{E}(t)\in \mathcal{L}(\mathcal{H}, \mathcal{H}')$, 
		where $\mathcal{H}'$ is the dual one for $\mathcal{H}$ and $\mathcal{H}'=X_{-\!1\!/\!2}\times H$,
		%with respect to the pivot space $X_{\frac{1}{2}}\times H$, 
		such that  $\mathcal{E}(t)$ is strictly positive and monotone increasing in $\mathcal{H}$, and 
		$$(-p_t^T, p^T)=\mathcal{E}(T-t)(w^T, w_t^T),$$
		where  $(w^T,  w^T_t)$ is the optimal state for system (\ref{001}) in the sense of (\ref{opt}), and $\mathcal{E}(\cdot)$ is the solution to the Riccati equation with initial condition $0$, that is,
		\begin{equation}\label{211}
		\left\{
		\begin{array}{l}
		\mathcal{E}_t=\mathcal{C}^*\mathcal{C}+(\mathcal{E}\mathcal{A}+\mathcal{A}^*\mathcal{E})-\mathcal{E}\mathcal{B}\mathcal{B}^*\mathcal{E},\quad in\; (0,+\infty),\\
		\mathcal{E}(0)=0,
		\end{array}
		\right.
		\end{equation}
in which $\mathcal{C}=[C, 0],\; \mathcal{B}=\left[\begin{array}{c} 0\\ B \end{array}\right]$.
		\end{lemma}
	
	{\bf Proof.}
	Since  $B\in \mathcal{L}(U, H)$ and $C\in \mathcal{L}(X_{1\!/\!2}, V)$ and $\mathcal{A}$ generates a $C_0$ semigroup on $\mathcal{H}$, by Theorem 2.1 (p. 393) in \cite{BDDM}, we obtain the 
unique existence of the Riccati operator  $\mathcal{E}(t)$.  In fact, it is a consequence of the fact that the optimality system (\ref{optimal system1}) has a unique solution, and the adjoint state $(-p_t, p)$ at time $t$ is a linear function of the state $(w, w_t)$ at time $t$.

Note that 
it can be checked directly that 
\begin{equation}\label{alge}
 \langle\mathcal{E}(T)(w_0, w_1),(w_0, w_1)\rangle_{\mathcal{H'}, \mathcal{H}}=\min_u J_0^T.
 \end{equation}
Thus, due to the (weak) observability of $(A, C)$ (see  (H2)), together with (\ref{alge}),  we get that
$$\langle \mathcal{E}(T)(w_0, w_1),(w_0, w_1)\rangle_{\mathcal{H'}, \mathcal{H}}\geq \frac{1}{c} \| (w_0, w_1)\|_{X_{{1}\!/\!{2}\!-\!{1}\!/\!{(2\eta)}}\times X_{-{1}\!/\!{(2\eta)}}}^2,$$ and hence
 $\mathcal{E}(\cdot)$ is strictly positive in $\mathcal{H}$.
 Moreover,  $\mathcal{E}(t)$ is monotone increasing. Indeed, for $t_1\leq t_2$, let $u_i,\; \varphi_i,\; i=1,2$ be the optimal control and trajectory in $[0, t_i]$. Then, we have 
 \begin{eqnarray*}
 \langle \mathcal{E}(t_2)(w_0, w_1), (w_0, w_1)\rangle_{\mathcal{H}', \mathcal{H}}&=&\min_u J_0^{t_2}=J_0^{t_2}(u_2)\geq J_0^{t_1}(u_2)\cr
 &&\geq \min_u J_0^{t_1}= J_0^{t_1}(u_1) =  \langle\mathcal{E}(t_1)(w_0, w_1), (w_0, w_1)\rangle_{\mathcal{H}', \mathcal{H}}.
 \end{eqnarray*}
Thus, $\mathcal{E}(t_2)\geq \mathcal{E}(t_1)$ for $t_2>t_1>0$. The proof is completed. \hfill$\Box$

Consider the corresponding  infinite-horizon quadratic cost functional associated with the control system (\ref{001}):
\begin{equation}\label{optinf}
\min J^\infty(u)=\frac{1}{2}\int_0^\infty [\|u(t)\|^2_U+\|Cw(t)\|_V^2]dt,\; u\in L^2(0,T; U).
\end{equation}

We have the following result on the well-posedness of the above infinite-horizon optimal control problem.
\begin{lemma}\label{l-2-3}
	 Assume that (H1) is fulfilled.  Then the set of admissible controls for problem (\ref{optinf}) is non-empty if the initial state $(w_0, w_1)\in \mathcal{D}(\mathcal{A}^k)=X_{(k+1)\!/\!2}\times X_{{k}\!/\!{2}}$ is sufficiently smooth satisfying $k\geq \frac{1}{\varrho}$. 
	\end{lemma}

{\bf Proof.~} 
%Based  on the weak observability of $(A, B^*)$,  we obtain that the set of admissible controls for the optimisation problem under consideration is non-empty if the initial state $(w_0, w_1)\in \mathcal{D}(\mathcal{A}^k)$ is sufficiently smooth satisfying $k\geq \frac{1}{\varrho}$. 
Due to the weak observability of $(A,B^*)$ as given  in (H1),   by the HUM method and  Proposition 3.25 in  \cite{dager}, p.43,  we obtain that for any given $(w_0,w_1)\in \mathcal{D}(\mathcal{A}^{\frac{1}{\varrho}})$,  there always exists a control $u(t)={\omega_{T_0}(t)}\in L^2((0,T_0], U),$ such that the solution to (\ref{001}) satisfying $(w(T_0), w_t(T_0))=0$ and  
\begin{equation}\label{221}
\int_0^{T_0} \|\omega_{T_0}(t)\|^2_Udt \leq c\|(w_0,w_1)\|_{\mathcal{D}(\mathcal{A}^{\frac{1}{\varrho}})}^2.
\end{equation}
Consider the control $\widetilde{u}$ on $ [0, \infty)$ that is equal to $\omega_{T_0}$ for $t\in(0, T_0]$ and is
identically zero for $t > T_0$. Thus,  by (\ref{221}), along with the well-posedness of control system (\ref{001}) and $C\in \mathcal{L}(X_{{1}\!/\!{2}}, H)$, we have
\begin{eqnarray}\label{112}
\min_u J^\infty(u)\leq J^\infty(\widetilde{u})
&=&\frac{1}{2}\int_0^{T_0} [\|\omega_{T_0}\|^2_U+\|C{w}(t)\|^2_V]dt\cr
&\leq& \widetilde{c} \|(w_0,w_1)\|_{\mathcal{D}(\mathcal{A}^{\frac{1}{\varrho}})}^2.
\end{eqnarray}

Hence, by \cite{BDDM}, we obtain the set of admissible controls for (\ref{optinf}) is non-empty for $(w_0, w_1)\in \mathcal{D}(\mathcal{A}^{\frac{1}{\varrho}})$.\hfill$\Box$

%\begin{remark}
%	In the above proof, we obtain the non-empty set of admissible controls by using the weak observability of $(A, B^*)$.  If using the stabilization result as given in (\ref{3}), then we can get a slightly weaker result that 
%	$$\min_u J^\infty(u)\leq J^\infty(-B^*w_t)\leq c \|(w_0,w_1)\|_{\mathcal{D}(\mathcal{A}^{k})}^2,\; k>\frac{1}{\varrho}. $$f
%	
%	\end{remark}

 Let 
 $(\hat{w}, \hat{w}_t)$ and $\hat{u}$ be the optimal state and control for system  (\ref{001})  respect to  (\ref{optinf}), and 
by  \cite{lixun},  we see that  $\hat{w}$ and $\hat{u}$ satisfies

\begin{equation}\label{osinf}
\left\{
\begin{array}{l}
\hat{w}_{tt} + A\hat{w} =   B\hat{u}, \\
w(0)=w_0, \quad w_t(0)=w_1,\\
\hat u= - B^*{\hat{p}}, \\
\hat{p}_{tt} + A\hat{p}=C^*C \hat{w}, \\
\hat{p}(t)\to 0,\; \hat{p}_t(t)\to 0,\quad t\to \infty.
\end{array}\right.
\end{equation}
\begin{proposition}\label{pro2.5} Assume that (H1) holds true. Then
there exists a unique minimal solution
	$$\hat{E}\in \mathcal{L}(\mathcal{D}(\mathcal{A}^{\frac{1}{\varrho}}), \mathcal{D}(\mathcal{A}^{\frac{1}{\varrho}})' )$$
	 of the  algebraic Riccati equation 
	 \begin{equation}\label{arica}
	 \mathcal{C}^*\mathcal{C}+(\hat{E}\mathcal{A}+\mathcal{A}^*\hat{E})-\hat{E}\mathcal{B}\mathcal{B}^*\hat{E}=0,
	 \end{equation}
	   such that 
	$(-\hat{p}_t, \hat{p})=\hat{E}(\hat{w}, \hat{w}_t).$
\end{proposition}

Thus, we can get the  Riccati-based optimal feedback control law for the infinite horizon problem. 
\begin{equation}\label{ricca2}
(0, \hat{u}(t))=-\mathcal{B}^*\hat{E} (\hat{w}, \hat{w}_t)=(0, -B^*\hat{p}(t)).
\end{equation}
and similarly we get that for  $ (w_0, w_1)\in \mathcal{D}(\mathcal{A}^{\frac{1}{\varrho}})$,
$$\quad \langle\hat{E}(w_0, w_1),(w_0, w_1)\rangle_{\mathcal{D}(\mathcal{A}^{\frac{1}{\varrho}})', \mathcal{D}(\mathcal{A}^{\frac{1}{\varrho}})}=\min_u J^\infty(u).$$

\subsection{Main results}
Under the assumptions (H1) and  (H2), in order to discuss the large time behaviour of system (\ref{001}) with Riccati-based optimal feedback control, we estimate the upper and lower bounds of $$\langle\hat{E} (w_0,w_1),(w_0,w_1)\rangle_{\mathcal{D}(\mathcal{A}^{\frac{1}{\varrho}})', \mathcal{D}(\mathcal{A}^{\frac{1}{\varrho}})},$$ and based on which, the explicit slow decay rate can be derived.

{\begin{theorem}\label{th1}
		Suppose that  (H1) and (H2) are satisfied. Then

		(1). There exists  constants $c_j>0,\; j=1,2$ such that
		
			\begin{eqnarray}\label{lowupper}	
		c_1{\| (w_0, w_1)\|_{X_{{1}\!/\!{2}-{1}\!/\!{(2\eta)}}\times X_{-{1}\!/\!{(2\eta)}}}^2}&\leq&  (\hat{E} (w_0,w_1),(w_0,w_1))_{\mathcal{D}(\!\mathcal{A}^{\frac{1}{\varrho}}\!)', \mathcal{D}(\!\mathcal{A}^{\frac{1}{\varrho}}\!)}\cr
		&\leq& c_2\|(w_0, w_1)\|^2_{\mathcal{D}(\mathcal{A}^\frac{1}{\varrho})}, 
		\end{eqnarray}
		where   $\varrho$ and $\eta$ is given as in  (H1) and  (H2), respectively.
		
%		\begin{equation}\label{lowupper}
%		c_1\|(w_0,w_1)\|^2_{\mathcal{H}}\leq  \langle\hat{E} (w_0,w_1),(w_0,w_1)\rangle_{\mathcal{D}(\mathcal{A}^{\frac{1}{\varrho}})', \mathcal{D}(\mathcal{A}^{\frac{1}{\varrho}})}\leq c_2\|(w_0, w_1)\|^2_{\mathcal{D}(\mathcal{A}^{\frac{1}{\varrho}})},
%		\end{equation}
%		where $c_j,\; j=1,2$ are some positive constants and $\varrho>0$ is given as in  (H1).
		
		(2).  There exists  a constant $M>0$ such that
		\begin{equation}\label{224}
	\|(\hat{w}(t), \hat{w}_t(t))\|_{\mathcal{H}}^2\leq 
	{M}{(t+1)^{-\frac{s\eta\varrho}{\varrho+\eta}}}\|(w_0,w_1)\|_{X_{({1}\!/\!{\varrho}+{1}\!/\!{\eta}+s+1)/2}\times X_{({1}\!/\!{\varrho}+{1}\!/\!{\eta}+s)/2}}^2,\; s>0,
	\end{equation}
		where $(\hat{w}(t), \hat{w}_t(t))$ is the solution to system (\ref{001}) with the Riccati-based optimal feedback control (\ref{ricca2}), $\varrho>0$ and $\eta>0$ is given as in (H1) and  (H2), respectively.

\end{theorem}

\begin{remark}
By the proof for Theorem \ref{th1}, we  find that the lower bound  in (\ref{lowupper}) is for $(w_0, w_1)\in X_{{1}\!/\!{2}-{1}\!/\!{(2\eta)}}\times X_{-{1}\!/\!{(2\eta)}}$ due to the weak observability of $(A, C)$ given as in (H2), while the upper bound is  for $(w_0, w_1)\in \mathcal{D}(\mathcal{A}^{\frac{1}{\varrho}})$ mainly caused by the weak observability of $(A,B^*)$ (see (H1)). Both of these two weak observability hypothesises determine the slow decay rate of the closed-loop system with the  Riccati-based optimal feedback control.
\end{remark}

%\begin{theorem}\label{th42}
%		Suppose that  (H1) and (H2) are satisfied. 
%	Let $(\hat{w}(t), \hat{w}_t(t))$ be the solution to system (\ref{001}) with the Riccati-based feedback control law, then there exist positive constant $M>0$ such that
%	$$ \| (\hat{w}(t), \hat{w}_t(t))\|_{\mathcal{H}}^2\leq {M}{(t+1)^{-\frac{s}{k}}} \|(w_0,w_1)\|_{\mathcal{D}(\mathcal{A}^{k+s})}^2,\quad k>\frac{1}{\varrho},\; s>0.$$
%\end{theorem}

\begin{remark}
		if assuming that the  exact  observability of $(A,B^*)$ and $(A,C)$  are fulfilled, i.e., there exist positive constants $T_0,\; T_1$ and $c$ such that 
	\begin{equation}\label{exact}
		c\int_0^{T_0}\|B^*w_t\|_{U'}^2dt\geq\| (w_0, w_1)\|^2_{\mathcal{H}},\quad 
			c\int_0^{T_1}\|C{w}\|_V^2dt\geq \| (w_0, w_1)\|^2_{\mathcal{H}},
		\end{equation}
		 we can see from (\ref{224}) that the system can be stabilized exponentially under Riccati-based optimal feedback control  for this case.  However, it should be noted that these strong assumptions only can be satisfied under suitable geometric conditions on the medium in which waves propagate. For instance, for waves on planar networks one needs a tree-like graph with control or observation in  all the free extreme edges except for one. Besides, for the multi-dimensional wave equation, the GCC has to be satisfied for control areas.
		
	If either  of the exact observation of $(A,B^*)$ and $(A,C)$    is fulfilled, then the exponential decay no longer holds, but the polynomial decay  can still be achieved. Specifically,  the  decay rates can be also derived from (\ref{224}) and given as $(t+1)^{-s\eta}$ for $(w_0,w_1)\in X_{({1}\!/\!{\eta}+s+1)/2}\times X_{({1}\!/\!{\eta}+s)/2}$ and $(t+1)^{-s{\varrho}}$ for $X_{({1}\!/\!{\varrho}+s+1)/2}\times X_{({1}\!/\!{\varrho}+s)/2}$, respectively. 
\end{remark}

\begin{remark}\label{remark1}
	If some other kinds of weak observability hypothesises are fulfilled, the corresponding slow decay rates can also be deduced similarly. For instance, If  (H1) and (H2) are replaced by the following much weaker ones:
	
$$
	c\int_0^T\|B^*w_t\|^2_{U'}dt\geq \| (w_0, w_1)\|_{\mathcal{D}(e^{-a\mathcal{A}})}^2$$
	and
	$$
	c\int_0^T\|Cw\|^2_{V}dt\geq \| (w_0, w_1)\|_{\mathcal{D}(e^{-b\mathcal{A}})}^2$$
	where 
	$$ \|(w_0, w_1)\|_{\mathcal{D}(e^{-\alpha\mathcal{A}})}^2:=\sum_{j\geq 1}\limits
e^{-2\alpha \lambda_n}(\lambda_n^2a_n^2+b_n^2),$$
then similar to the proof for Theorem \ref{th1},
the lower and  upper bound become
\begin{equation}\label{exp1}
 c_1\| (w_0, w_1)\|_{\mathcal{D}(e^{-b\mathcal{A}})}^2\leq  \langle\hat{E} (w_0,w_1),(w_0,w_1)\rangle_{\mathcal{D}(e^{k\mathcal{A}})', \mathcal{D}(e^{k\mathcal{A}})}\!\leq\! c_2\|(w_0, w_1)\|^2_{\mathcal{D}(e^{a\mathcal{A}})}.
 \end{equation}
 
Moreover, we can obtain the following slow decay rate for the system with Riccati-based optimal feedback control (\ref{ricca2}), that is,
%Thus, similar to the proof for Theorem \ref{th1}, we can get that
	
	\begin{equation}\label{exp2}
	 \| (\hat{w}(t), \hat{w}_t(t))\|_{\mathcal{H}}^2\leq {M}{(t+1)^{-\frac{s}{a+b}}} \|(w_0,w_1)\|_{\mathcal{D}(e^{(a+b+s)\mathcal{A}})}^2,\; s>0.
	 \end{equation}
%	
%	It should be noted that although the logarithmic decay of the closed-loop system can also be obtained when the initial data belongs to $\mathcal{D}(\mathcal{A})$, it is not sufficient to guarantee that the  set of admissible controls is non-empty. Indeed, the logarithmic decay rate is too slow to make $\min_u\limits J^\infty (u)$ be bounded.
	\end{remark}

Let $w^T(t)$, $u^T(t)$ be the optimal solution and control to the following problem 
\begin{equation}\label{opt-z}
\min \widetilde{J}_0^T(u)=\frac{1}{2}\int_0^T [\|u(t)\|^2_U
+\|Cw(t)-z\|_V^2]dt,\; u\in L^2(t,T; U)
\end{equation}
where $w$ and $u$ satisfy equation (\ref{001}).

Let  $\bar{w}, \bar{u}$ be the optimal solution and control to the steady optimal control problem 
\begin{equation}\label{staopt}
\min_u \widetilde{J}_s=\|u\|^2_U+\|Cw-z\|_V^2,
\end{equation}
subject to $A w=Bu$.  
For the stationary optimal control problem (\ref{staopt}), $(\bar u, \bar w)$ satisfy
$$A \bar{w}=B \bar{u},\; \langle\bar{u},  v\rangle+\langle C\bar{w}-z, C\varphi\rangle=0,\; {\rm for\; every\;} v, \varphi: A\varphi=Bv.$$

Thus, $C^*(C\bar{w}-z)\in Ker(A)^\perp$ and hence there exists some $\bar{p}$ satisfying 
$$A^* \bar{p}=A\bar p=C^*(C\bar{w}-z).$$

Hence, by the above and (\ref{optimal system1}), we obtain the following OS system
\begin{equation}\label{403}
\left\{
\begin{array}{l}
(w^T-\bar{w})_{tt}+A(w^T-\bar{w})=B(u^T-\bar{u}),\\
(p^T-\bar{p})_{tt}+ A(p^T-\bar{p})=C^*C (w^T-\bar{w}), \\
u^T-\bar{u}= - B^*(p^T-\bar{p}), \\
w^T(0)=w_0, \quad w_t^T(0)=w_1,\\
p^T(T)=p_t^T(T)=0.
\end{array}
\right.
\end{equation}

Let us discuss the turnpike property of the optimal control problem, that is, identifying that to which extent do  $u^T(t), w^T(t)$  approximate the stationary ones $\bar{u}, \bar{w}$ as $T\to \infty$.  In fact, we have the following result.

\begin{theorem}\label{t-4-3}
	Assume that (H1) and (H2) hold true. Then for any $ (w_0-\bar{w}, w_1) \in\mathcal{D}(\mathcal{A}^k),\; k>0$ and  $\bar{p}\in {X_{({k}+1)\!/\!{2}}},\; {k}\geq 1$, it holds that as $T\to\infty$,
	\begin{equation}\label{281}
	\frac{1}{T}\int_0^T(\|C(w^T-\bar{w})\|_{H}^2+\|(u^T-\bar{u})\|_H^2)dt\to 0,
	\end{equation}
and
	\begin{equation}\label{282}
	\|\frac{1}{T}\int_0^T (w^T(t)-\bar{w})dt\|_{X_{{1}\!/\!{2}}}^2 \to 0.
	\end{equation}
\end{theorem}

\begin{remark}
We see that 	when choosing  (\ref{exact}) instead of  (H1) and (H2), i.e., the exact observability of $(A,B^*)$ and $(A,C)$, the result in Theorem \ref{t-4-3} is consistent with the  averaged turnpike property as given in \cite{zuazua2}. Thus, in terms of averaged turnpike property, the result in \cite{zuazua2} can be considered as a special case of Theorem \ref{t-4-3} given above.
	\end{remark}
\begin{remark}
	In the above theorem, the averaged turnpike property is obtained under sufficiently smooth initial states $(w_0-\bar{w}, w_1)$ and $\bar{p}$. Compared to the exponential-type point-wise turnpike property under exact observability of $(A,B^*)$ and $(A,C)$, it is reasonable to expect the  polynomial point-wise turnpike property under the weak observability hypothesises (H1) and (H2). However, the point-wise turnpike property is still an open problem. In fact, the proof based on energy estimates along with the properties of  Riccati operator   proposed in \cite{zuazua2}  is difficult to apply to the case with weak observability as given in (H1) and (H2) under consideration. The main difficulty in the energy estimates is  caused by the slow decay rate  which is determined by the weak observability. As we know, the slow decay rate is non-uniform and  always dependent on the regularity of the initial states.  This is different from the case in \cite{zuazua2}, where the exact observability is  fulfilled and the uniform exponential decay rates always hold for any initial data in state space.
	 
\end{remark}

\section{Proof of main results}\label{333}
\subsection{Proof of Proposition \ref{pro2.5}}

	Following the proof in Lemma  \ref{l-2-3}, due to (H1),  for  $(w_0, w_1)\in \mathcal{D}(\mathcal{A}^{\frac{1}{\varrho}})$, we always can find a control  $\widetilde{u}=\left\{
	\begin{array}{l}
	\omega_{T_0}(t),\quad  t\in (0, T_0], \\
	0,\quad t\in (T_0, T],
	\end{array}
	\right.$
	such that 
	$$(w(t), w_t(t))=0,\; t\geq T_0$$ and  $$\int_0^T\|\widetilde{u}\|^2_Udt =\int_0^{T_0} \|\omega_{T_0}(t)\|^2_Udt \leq c\|(w_0,w_1)\|_{\mathcal{D}(\mathcal{A}^{\frac{1}{\varrho}})}^2.$$
	Hence,  we  get that for $(w_0, w_1)\in \mathcal{D}(\mathcal{A}^{\frac{1}{\varrho}})$, 
	\begin{eqnarray}\label{uppb}
	\langle\mathcal{E}(T)(w_0,w_1), (w_0,w_1) \rangle_{\mathcal{D}(\mathcal{A}^{\frac{1}{\varrho}})', \mathcal{D}(\mathcal{A}^{\frac{1}{\varrho}})}&=&\min_u J_0^T\cr
	&\leq& 
	J_0^T(\widetilde{u})=\frac{1}{2}\int_0^{T_0} [\|\omega_{T_0}\|^2_U+\|C{w}(t)\|^2_V]dt\cr
	&\leq&  \widetilde{c} \|(w_0,w_1)\|_{\mathcal{D}(\mathcal{A}^{\frac{1}{\varrho}})}^2, \quad {\rm for \; all\,} T>T_0.
	\end{eqnarray}
	
	Note that for $(w_0,w_1), (z_0, z_1)\in \mathcal{D}(\mathcal{A}^{\frac{1}{\varrho}})$,  we have 
	\begin{eqnarray*}
		&&2{\rm Re }\langle\mathcal{E}(T) (w_0, w_1), (z_0, z_1)\rangle_{\mathcal{D}(\mathcal{A}^{\frac{1}{\varrho}})', \mathcal{D}(\mathcal{A}^{\frac{1}{\varrho}})}\cr
		&=&\langle\mathcal{E}(T) (w_0+z_0, w_1+z_1), (w_0+z_0, w_1+z_1)\rangle_{\mathcal{D}(\mathcal{A}^{\frac{1}{\varrho}})', \mathcal{D}(\mathcal{A}^{\frac{1}{\varrho}})}\cr
		&&-\langle\mathcal{E}(T) (w_0, w_1), (w_0, w_1)\rangle_{\mathcal{D}(\mathcal{A}^{\frac{1}{\varrho}})', \mathcal{D}(\mathcal{A}^{\frac{1}{\varrho}})}\cr
		&&-\langle\mathcal{E}(T) (z_0, z_1), (z_0, z_1)\rangle_{\mathcal{D}(\mathcal{A}^{\frac{1}{\varrho}})', \mathcal{D}(\mathcal{A}^{\frac{1}{\varrho}})},
	\end{eqnarray*}
	\begin{eqnarray*}
		&&2{\rm Im }\langle\mathcal{E}(T) (w_0, w_1), (z_0, z_1)\rangle_{\mathcal{D}(\mathcal{A}^{\frac{1}{\varrho}})', \mathcal{D}(\mathcal{A}^{\frac{1}{\varrho}})}\cr
		&=&i\langle\mathcal{E}(T) (w_0+iz_0, w_1+iz_1), (w_0+iz_0, w_1+iz_1)\rangle_{\mathcal{D}(\mathcal{A}^{\frac{1}{\varrho}})', \mathcal{D}(\mathcal{A}^{\frac{1}{\varrho}})}\cr
		&&-\langle\mathcal{E}(T) (w_0, w_1), (w_0, w_1)\rangle_{\mathcal{D}(\mathcal{A}^{\frac{1}{\varrho}})', \mathcal{D}(\mathcal{A}^{\frac{1}{\varrho}})}\cr
		&&-\langle\mathcal{E}(T) (iz_0, iz_1), (iz_0, iz_1)\rangle_{\mathcal{D}(\mathcal{A}^{\frac{1}{\varrho}})', \mathcal{D}(\mathcal{A}^{\frac{1}{\varrho}})},
	\end{eqnarray*}
	which along with the monotone increasing property of $\mathcal{E}(t)$ and (\ref{uppb}), yields that  the limit $$\lim_{T\to\infty}\limits\langle\mathcal{E}(T)(w_0, w_1), (z_0, z_1)\rangle_{\mathcal{D}(\mathcal{A}^{\frac{1}{\varrho}})', \mathcal{D}(\mathcal{A}^{\frac{1}{\varrho}})}$$ exists for all $ (w_0, w_1), (z_0, z_1)\in\mathcal{D}(\mathcal{A}^{\frac{1}{\varrho}})$. 
	By the Uniform Boundedness Theorem, it follows that $\mathcal{E}(\cdot)$ is bounded in $\mathcal{D}(\mathcal{A}^{\frac{1}{\varrho}})$. 
	%Set $$\gamma(x)=\lim_{T\to\infty}(\mathcal{E}(T) (w_0, w_1), (w_0, w_1))_{\mathcal{D}(\mathcal{A}^{\frac{1}{\varrho}})', \mathcal{D}(\mathcal{A}^{\frac{1}{\varrho}})},\quad {\rm for\; all\;} (w_0, w_1)\in \mathcal{D}(\mathcal{A}^{\frac{1}{\varrho}}).$$
	Thus, for $ (w_0, w_1), (z_0, z_1)\in\mathcal{D}(\mathcal{A}^{\frac{1}{\varrho}})$, we define 
	$$\lim_{T\to\infty}\langle\mathcal{E}(T)(w_0, w_1), (z_0, z_1)\rangle_{\mathcal{D}(\mathcal{A}^{\frac{1}{\varrho}})', \mathcal{D}(\mathcal{A}^{\frac{1}{\varrho}})}:=\langle\hat{E}(w_0, w_1), (z_0, z_1)\rangle_{\mathcal{D}(\mathcal{A}^{\frac{1}{\varrho}})', \mathcal{D}(\mathcal{A}^{\frac{1}{\varrho}})}.
	$$
	It follows that for $
	(w_0, w_1)\in\mathcal{D}(\mathcal{A}^{\frac{1}{\varrho}})$,
	$$\lim_{T\to\infty}\langle\mathcal{E}(T)(w_0, w_1), (w_0, w_1)\rangle_{\mathcal{D}(\mathcal{A}^{\frac{1}{\varrho}})', \mathcal{D}(\mathcal{A}^{\frac{1}{\varrho}})}=\langle\hat{E}(w_0, w_1), (w_0, w_1)\rangle_{\mathcal{D}(\mathcal{A}^{\frac{1}{\varrho}})', \mathcal{D}(\mathcal{A}^{\frac{1}{\varrho}})}.$$
	
	Hence,  for $ (w_0, w_1)\in \mathcal{D}(\mathcal{A}^{\frac{1}{\varrho}})$, the limit $ \lim_{T\to\infty}\limits \mathcal{E}(T)(w_0, w_1)=\hat{E}(w_0, w_1)$ exists in $\mathcal{D}(\mathcal{A}^{\frac{1}{\varrho}})'$.
	Then following the proof of Proposition 2.2  in \cite{BDDM} (p.483), we get that $\hat{E}$ satisfies 
	(\ref{arica}) and $\hat{E}\leq X$ for any solution $X$ of (\ref{arica}) and hence $\hat{E}$ is the minimal solution to the algebraic Riccati equation (\ref{arica}). The proof is completed.\hfill$\Box$
\subsection{Slow decay rate (Proof of Theorem \ref{th1})}\label{proofall}

{First,  by (\ref{112}), we can easily obtain the upper bound as follows.

\begin{equation}\label{up300}
\langle\hat{E}(w_0, w_1),(w_0, w_1)\rangle_{\mathcal{D}(\mathcal{A}^{\frac{1}{\varrho}})', \mathcal{D}(\mathcal{A}^{\frac{1}{\varrho}})}=\min_u J^\infty\leq c_1 \|(w_0, w_1)\|^2_{\mathcal{D}(\mathcal{A}^{\frac{1}{\varrho}})}.
\end{equation}

}

In order to get the lower bound, the proof is given by the following two steps:

{\it Step 1.} Divide $w^T$ in system (\ref{optimal system1}) by $w^T=y+z$, where
$y$ satisfies
\begin{equation}\label{ysystem}
\left\{
\begin{array}{l}
y_{tt} + Ay =   Bu^T, \quad 0<t< T, \\
y(0)=0, \quad y_t(0)=0,\\
u^T= - B^*p^T,
\end{array}\right.
\end{equation}
and $z$ satisfies
\begin{equation}\label{zsystem}
\left\{
\begin{array}{l}
z_{tt} + Az = 0,\quad 0<t< T, \\
z(0)=w_0, \quad z_t(0)=w_1.
\end{array}\right.
\end{equation}
{By the hypothesis (H2), we get  directly that  there exists a constant $c>0$ satisfying
\begin{equation}\label{03+}
{\| (w_0, w_1)\|_{X_{{1}\!/\!{2}-{1}\!/\!{(2\eta)}}\times X_{-{1}\!/\!{(2\eta)}}}^2}\leq c\int_0^T\|Cz\|_V^2dt,
\end{equation}
where $z$ satisfies (\ref{zsystem}).}

We also have the following estimate for (\ref{ysystem}).
\begin{equation}\label{110}
\|(y,y_t)\|_{\mathcal{ H}}^2 \leq e^{t} \int_0^t \|Bu^T(s)\|^2_{H}ds,\quad 0<t\leq T.
\end{equation}
In fact, set the energy function $E(t)=\frac{1}{2}(\|y\|_{X_{{1}\!/\!{2}}}^2+ \|y_t\|_H^2)$.
Then, differentiating $E(t)$ by $t$, together with (\ref{ysystem}), we get 
$$\frac{d E(t)}{dt}= \langle y_t, Bu^T\rangle.$$

Hence,
\begin{eqnarray}
E(t)&\leq& \frac{1}{2}\int_0^t(\|Bu^T\|^2_{H}+\|y_t\|_H^2)ds\cr
&\leq&  \frac{1}{2}\int_0^t\|Bu^T(s)\|^2_{H}ds+ \int_0^t E(s)ds,\;\; 0<t\leq T.
\end{eqnarray}

Using Gronwall's inequality, we obtain
\begin{equation}
E(t)\leq \frac{1}{2}e^{t}\int_0^t\|Bu^T(s)\|_H^2ds,\quad 0<t\leq T,
\end{equation}
which leads to (\ref{110}).

{\it Step 2.\quad  we will show that 
	$$\|(w_0,w_1)\|_{X_{{1}\!/\!{2}-{1}\!/\!{(2\eta)}}\times X_{-{1}\!/\!{(2\eta)}}}^2\leq {c} \langle\hat{E} (w_0,w_1),(w_0,w_1)\rangle_{\mathcal{D}(\mathcal{A}^{\frac{1}{\varrho}})', \mathcal{D}(\mathcal{A}^{\frac{1}{\varrho}})},$$
	where $c>0$ is some constant.
}

 Note that $w^T=y+z$.
Thus, 
\begin{eqnarray}\label{308-}
{\| (w_0, w_1)\|_{X_{{1}\!/\!{2}-{1}\!/\!{(2\eta)}}\times X_{-{1}\!/\!{(2\eta)}}}^2}&\leq& c\int_0^T\|Cz\|_V^2dt\cr
&\leq& c\int_0^T(\|Cw^T\|_V^2+\|Cy\|_V^2)dt\cr
&\leq&  c\int_0^T(\|Cw^T\|^2_V+\|y\|_{X_{{1}\!/\!{2}}}^2)dt\cr
&\leq &c\int_0^T(\|Cw^T\|_V^2+ e^t \int_0^t\|Bu^T(s)\|_H^2ds)dt\cr
&\leq&c\cdot \max \{1, e^T-1\}\int_0^T (\|Cw^T\|_V^2+ \|Bu^T\|_H^2)dt\cr
&\leq&c_T \int_0^T (\|Cw^T\|_V^2+ \|u^T\|_U^2)dt\cr
&=& c_T \langle\mathcal{E}(T)(w_0, w_1),(w_0, w_1)\rangle_{\mathcal{H}', \mathcal{H}},
\end{eqnarray}
where $c_T=c\cdot \max \{1, e^T-1\}$.
Thus, by (\ref{308-}), we have
\begin{equation}\label{308}
\langle\mathcal{E}(T)(w_0, w_1),(w_0, w_1)\rangle_{\mathcal{H}', \mathcal{H}}\geq   \frac{1}{c_T}{\| (w_0, w_1)\|_{X_{{1}\!/\!{2}-{1}\!/\!{(2\eta)}}\times X_{-{1}\!/\!{(2\eta)}}}^2}. 
\end{equation}

Note that by Proposition \ref{pro2.5}, we know that $\mathcal{E}(T)$ is monotone increasing and bounded in $\mathcal{D}(\mathcal{A}^{\frac{1}{\varrho}})$ and $ \lim_{T\to\infty}\limits \mathcal{E}(T)(w_0, w_1)=\hat{E}(w_0, w_1)$ exists in $\mathcal{D}(\mathcal{A}^{\frac{1}{\varrho}})'$. Hence,  for $(w_0,w_1)\in \mathcal{D}(\mathcal{A}^{\frac{1}{\varrho}})$,
\begin{eqnarray}\label{low309}
&&\langle\hat{E}(w_0, w_1),(w_0, w_1)\rangle_{\mathcal{D}(\mathcal{A}^{\frac{1}{\varrho}})', \mathcal{D}(\mathcal{A}^{\frac{1}{\varrho}})}\cr
&\geq&  \langle\mathcal{E}(T)(w_0, w_1),(w_0, w_1)\rangle_{\mathcal{D}(\mathcal{A}^{\frac{1}{\varrho}})', \mathcal{D}(\mathcal{A}^{\frac{1}{\varrho}})}\cr
&\geq& 
\frac{1}{c_T}{\| (w_0, w_1)\|_{X_{{1}\!/\!{2}-{1}\!/\!{(2\eta)}}\times X_{-{1}\!/\!{(2\eta)}}}^2},\quad  \forall T>0.
\end{eqnarray}

Now, we consider the decay rate of the closed-loop system with Riccati-based optimal feedback control (\ref{ricca2}). A direct calculation yields 
\begin{eqnarray}\label{309}
\frac{d \langle\hat{E} (\hat{w}, \hat{w}_t),(\hat{w}, \hat{w}_t)\rangle_{\mathcal{D}(\mathcal{A}^{\frac{1}{\varrho}})', \mathcal{D}(\mathcal{A}^{\frac{1}{\varrho}})}}{dt}
	=
	-\|B^*\hat{p}\|_{U'}^2-\|C\hat{w}\|_V^2. 
\end{eqnarray}

Integrating the above  from $0$ to $T$, we have 
\begin{eqnarray}\label{402}
&&\langle\hat{E} (\hat{w}(T), \hat{w}_t(T)),(\hat{w}(T), \hat{w}_t(T))\rangle_{\mathcal{D}(\mathcal{A}^{\frac{1}{\varrho}})', \mathcal{D}(\mathcal{A}^{\frac{1}{\varrho}})}-\langle\hat{E} ({w}_0, {w}_1),({w}_0, {w}_1)\rangle_{\mathcal{D}(\mathcal{A}^{\frac{1}{\varrho}})', \mathcal{D}(\mathcal{A}^{\frac{1}{\varrho}})}\cr
&=&- \int_0^T(\|B^*\hat{p}\|_{U'}^2+\|C\hat{w}\|_V^2)dt.
\end{eqnarray}

Note from (\ref{308}) that 
\begin{eqnarray}
&&\int_0^T (\|B^*\hat{p}\|_{U'}^2+\|C\hat{w}\|_V^2)dt\geq  \langle\mathcal{E}(T)(w_0, w_1),(w_0, w_1)\rangle\cr
 &\geq& \frac{1}{c_T}{\| (w_0, w_1)\|_{X_{{1}\!/\!{2}-{1}\!/\!{(2\eta)}}\times X_{-{1}\!/\!{(2\eta)}}}^2}.
\end{eqnarray}

Thus,
\begin{eqnarray}\label{125+}
&&\langle\hat{E} (\hat{w}(T), \hat{w}_t(T)),(\hat{w}(T), \hat{w}_t(T))\rangle_{\mathcal{D}(\mathcal{A}^{\frac{1}{\varrho}})', \mathcal{D}(\mathcal{A}^{\frac{1}{\varrho}})}-\langle\hat{E} ({w}_0, {w}_1),({w}_0, {w}_1)\rangle_{\mathcal{D}(\mathcal{A}^{\frac{1}{\varrho}})', \mathcal{D}(\mathcal{A}^{\frac{1}{\varrho}})}\cr
&\leq& -\frac{1}{c_T}{\| (w_0, w_1)\|_{X_{{1}\!/\!{2}-{1}\!/\!{(2\eta)}}\times X_{-{1}\!/\!{(2\eta)}}}^2}
\end{eqnarray}	
	where  $\varrho>0$ and $\eta>0$ is given as   (H1) and  (H2), respectively.
	
	Then by interpolation, we have
	$$ \|(w_0,w_1)\|_{\mathcal{D}(\mathcal{A}^{\frac{1}{\varrho}})}^2\leq \| (w_0, w_1)\|_{X_{{1}\!/\!{2}-{1}\!/\!{(2\eta)}}\times X_{-{1}\!/\!{(2\eta)}}}^{\frac{2s\eta}{1+\frac{1}{\varrho}\eta+s\eta}}\|(w_0,w_1)\|_{\mathcal{D}(\mathcal{A}^{\frac{1}{\varrho}+s})}^{\frac{2(1+\frac{1}{\varrho}\eta)}{1+\frac{1}{\varrho}\eta+s\eta}},\; s>0,$$
	and hence
	$${\| (w_0, w_1)\|_{X_{{1}\!/\!{2}-{1}\!/\!{(2\eta)}}\times X_{-{1}\!/\!{(2\eta)}}}^2}\geq \frac{\|(w_0,w_1)\|_{\mathcal{D}(\mathcal{A}^{\frac{1}{\varrho}})}^{\frac{2(1+\frac{1}{\varrho}\eta+s\eta)}{s\eta}}}{\|(w_0,w_1)\|_{\mathcal{D}(\mathcal{A}^{\frac{1}{\varrho}+s})}^{\frac{2(1+\frac{1}{\varrho}\eta)}{s\eta}} },\; s>0.$$
	Thus, by (\ref{125+}), we have
	
	\begin{eqnarray}
	&&\langle\hat{E} (\hat{w}(T), \hat{w}_t(T)),(\hat{w}(T), \hat{w}_t(T))\rangle_{\mathcal{H}', \mathcal{H}}-\langle\hat{E} ({w}_0, {w}_1),({w}_0, {w}_1)\rangle_{\mathcal{H}', \mathcal{H}}\cr&\leq& -\frac{1}{c_T} \frac{\|(w_0,w_1)\|_{\mathcal{D}(\mathcal{A}^{\frac{1}{\varrho}})}^{\frac{2(1+\frac{1}{\varrho}\eta+s\eta)}{s\eta}}}{\|(w_0,w_1)\|_{\mathcal{D}(\mathcal{A}^{\frac{1}{\varrho}+s})}^{\frac{2(1+\frac{1}{\varrho}\eta)}{s\eta}} },\; s>0.
	\end{eqnarray}
	Therefore,
	\begin{eqnarray}
	&&\frac{\langle\hat{E} (\hat{w}(T), \hat{w}_t(T)),(\hat{w}(T), \hat{w}_t(T))\rangle_{\mathcal{H}', \mathcal{H}}}{\|(w_0,w_1)\|_{\mathcal{D}(\mathcal{A}^{\frac{1}{\varrho}+s})}^2}-\frac{\langle\hat{E} ({w}_0, {w}_1),({w}_0, {w}_1)\rangle_{\mathcal{H}', \mathcal{H}}}{\|(w_0,w_1)\|_{\mathcal{D}(\mathcal{A}^{\frac{1}{\varrho}+s})}^2}\cr
	&\leq& -\frac{1}{c_T} \left[\frac{\|(w_0,w_1)\|_{\mathcal{D}(\mathcal{A}^{\frac{1}{\varrho}})}^2}{\|(w_0,w_1)\|_{\mathcal{D}(\mathcal{A}^{\frac{1}{\varrho}+s})}^2}\right]^{\frac{1+\frac{1}{\varrho}\eta+s\eta}{s\eta}},\; s>0.
	\end{eqnarray}

	Thus, by (\ref{up300}) and the monotone decreasing property of $\langle\hat{E} ({w}, {w}_t),({w}, {w}_t)\rangle$, we have
	\begin{eqnarray}\label{128+++}
	&&\frac{\langle\hat{E} (\hat{w}(T), \hat{w}_t(T)),(\hat{w}(T), \hat{w}_t(T))\rangle_{\mathcal{H}', \mathcal{H}}}{\|(w_0,w_1)\|_{\mathcal{D}(\mathcal{A}^{\frac{1}{\varrho}+s})}^2}-\frac{\langle\hat{E} ({w}_0, {w}_1),({w}_0, {w}_1)\rangle_{\mathcal{H}', \mathcal{H}}}{\|(w_0,w_1)\|_{\mathcal{D}(\mathcal{A}^{\frac{1}{\varrho}+s})}^2}\cr
	&\leq& -\frac{1}{c_2^2c_T} [\frac{\langle\hat{E} (\hat{w}(T), \hat{w}_t(T)),(\hat{w}(T), \hat{w}_t(T))\rangle_{\mathcal{H}', \mathcal{H}}}{\|(w_0,w_1)\|_{\mathcal{D}(\mathcal{A}^{\frac{1}{\varrho}+s})}^2}]^{\frac{1+\frac{1}{\varrho}\eta+s\eta}{s\eta}},\;  s>0.
	\end{eqnarray}
	
In order to obtain the explicit decay rate of closed-loop system, let us introduce the following result in Ammari and Tucsnak \cite{ammari-tucsnak}.
\begin{lemma}\label{lemma4-1}
	Let $\{a_m\}_{m=1}^\infty$ be a sequence of positive number
	satisfying
	\begin{equation}
	a_{m+1}\leq a_m-C(a_{m+1})^{2+\alpha},\quad \forall m\geq 1,
	\end{equation}
	for some constants $C>0$ and $\alpha>-1$. Then there exists a
	positive constant $M_{C,\alpha}$ such that $$a_m\leq
	\frac{M_{C,\alpha}}{(m+1)^{\frac{1}{1+\alpha}}}.$$
\end{lemma}

By the above Lemma, together with (\ref{128+++}), it is easy to get
that 
$$
\frac{\langle\hat{E} (\hat{w}(t), \hat{w}_t(t)),(\hat{w}(t), \hat{w}_t(t))\rangle_{\mathcal{D}(\mathcal{A}^{\frac{1}{\varrho}})', \mathcal{D}(\mathcal{A}^{\frac{1}{\varrho}})}}{\|(w_0,w_1)\|_{\mathcal{D}(\mathcal{A}^{{\frac{1}{\varrho}}+s})}^2}\leq {M}{(t+1)^{-\frac{s\eta}{1+\frac{1}{\varrho}\eta}}},\; s>0.
$$

Note that by (\ref{low309}), we know
$$ \langle\hat{E} (\hat{w}(t), \hat{w}_t(t)),(\hat{w}(t), \hat{w}_t(t))\rangle_{\mathcal{D}(\mathcal{A}^{\frac{1}{\varrho}})', \mathcal{D}(\mathcal{A}^{\frac{1}{\varrho}})}\geq c_1  \| (\hat{w}(t), \hat{w}_t(t))\|_{X_{{1}\!/\!{2}-{1}\!/\!{(2\eta)}}\times X_{-{1}\!/\!{(2\eta)}}}^2.$$

So,
	$$
	\|(\hat{w}(t), \hat{w}_t(t))\|_{X_{{1}\!/\!{2}-{1}\!/\!{(2\eta)}}\times X_{-{1}\!/\!{(2\eta)}}}^2\leq 
	{M}{(t+1)^{-\frac{s\eta}{1+\frac{1}{\varrho}\eta}}}\|(w_0,w_1)\|_{\mathcal{D}(\mathcal{A}^{\frac{1}{\varrho}+s})}^2,\;  s>0,
	$$
	and hence
	\begin{equation}
	\|(\hat{w}(t), \hat{w}_t(t))\|_{\mathcal{H}}^2\leq 
	{M}{(t+1)^{-\frac{s\eta\varrho}{\varrho+\eta}}}\|(w_0,w_1)\|_{X_{({1}\!/\!{\varrho}+{1}\!/\!{\eta}+{s}+1)/2}\times X_{({1}\!/\!{\varrho}+{1}\!/\!{\eta}+{s})/2}}^2,\; s>0,
	\end{equation}
	where $\varrho>0$ and $\eta>0$ is given as in (H1) and  (H2), respectively.
	The proof is completed. \hfill$\Box$

\subsection{Averaged turnpike property (Proof of Theorem \ref{t-4-3})}\label{444}
Let us consider the  averaged turnpike property for the LQ optimal  control problems with (H1) and (H2) are fulfilled.
%For simplicity, we first choose $z=0$. Then it is obvious that $\bar{w}=\bar{u}=0$.
%We have the following result.
% Firstly,  based on Theorem \ref{th1}, we have the following result on the estimate for  $\mathcal{E}(T)-\hat{E}$.
Firstly, based on the weak observability of $(A,B^*)$, we obtain the following estimate.

\begin{lemma}\label{L-1.1}
	Assume that (H1) is fulfilled. Then there exists some positive constant $	\widetilde{c}$ such that 
	\begin{equation}\label{401}
	\|(p(0), p_t(0))\|_{\mathcal{D}(\mathcal{A}^{k})'}^2 
	\leq
	\widetilde{c}g_1(T)[\int_0^T\|B^*p\|_{H}^2dt+\int_0^T\|f\|_{X_{-{1}\!/\!{2}}}^2dt
	+\|p_{0T}\|_{H}^2 + \|p_{1T}\|_{X_{-{1}\!/\!{2}}}^2],
	\end{equation}
 where $k>0$, $g_1(T)= \left\{
 \begin{array}{ll} \frac{(T+1)^{-k\varrho+1}-1}{-k\varrho+1}& k\varrho\neq1,\\
 \ln(T+1)& k\varrho=1,
 \end{array}
 \right.$ and $p$ is any solution to the following inhomogeneous system
	\begin{equation}\label{109}
	\left\{
	\begin{array}{l}
	p_{tt}+A p=f,\\
	p(T)=p_{0T},\quad p_t(T)=p_{1T}.
	\end{array}
	\right.
	\end{equation}
\end{lemma}
{\bf Proof.~}By duality we get
\begin{equation}\label{406}
\langle p_t,y\rangle\mathop{|}\limits_0^T -  \langle p,y_t\rangle\mathop{|}\limits_0^T
+
\int_0^T \langle p,  y_{tt}+Ay \rangle \,dt=
\int_0^T \langle f,y \rangle\, dt.
\end{equation}
Let $y$ be the solution to equation (\ref{001}) with (\ref{002}). Using H\"{o}lder's inequality, we get
\begin{eqnarray*}
| \int_0^T \langle p,  y_{tt}+Ay \rangle \,dt| 
&=&| \int_0^T \langle p,  -BB^*y_t \rangle \,dt| \cr
&\leq& 
c(\int_0^T\|B^*p\|^2_{H}dt)^{\frac{1}{2}}(\int_0^T \|B^* y_t\|^2_{H})^{\frac{1}{2}}.
\end{eqnarray*}
Thus, by Lemma \ref{le1-11}, along with (\ref{20002+}),  we get that there exists some constant $c>0$ such that
\begin{eqnarray}\label{407}
&&| \int_0^T \langle p,  y_{tt}+Ay \rangle \,dt|\cr
&\leq& c (\int_0^T \| B^*p\|_{H}^2\, dt)^{\frac12} ( \frac{1}{2}\|(y_0, y_1)\|^2_{\mathcal{H}}+\frac{1}{2}c_1 (T+1)^{-{k} {\varrho} }\|(y_0, y_1)\|_{\mathcal{D}(\mathcal{A}^{{k}})}^2)^{\frac{1}{2}}\cr
&\leq&{c}  (\frac{1}{2}+\frac{1}{2}c_1 (T+1)^{-{k} {\varrho} })^{\frac{1}{2}} (\int_0^T \| B^*p\|_{H}^2\, dt)^{\frac12} \|(y_0, y_1)\|_{\mathcal{D}(\mathcal{A}^k)}.
\end{eqnarray}

Besides, we have 
\begin{eqnarray}\label{408}
\int_0^T \langle f,y\rangle \, dt
&\leq&  (\int_0^T \|f\|_{X_{-{1}\!/\!{2}}}^2)^{\frac{1}{2}} (\int_0^T \|y\|_{X_{{1}\!/\!{2}}}^2)^{\frac{1}{2}}\cr
&\leq& c (\int_0^T \|f\|_{X_{-{1}\!/\!{2}}}^2)^{\frac{1}{2}} (\int _0^T(t+1)^{-{k\varrho}} \|(y_0, y_1)\|_{\mathcal{D}(\mathcal{A}^k)}^2)^\frac{1}{2}\cr
&\leq&{c} (g_1(T))^{\frac{1}{2}} (\int_0^T \|f\|_{X_{-{1}\!/\!{2}}}^2)^{\frac{1}{2}}  \|(y_0, y_1)\|_{\mathcal{D}(\mathcal{A}^k)},\quad k>0.
\end{eqnarray}

Hence, by (\ref{406}), (\ref{407}) and (\ref{408}), along with Lemma \ref{le1-11},  we obtain
\begin{eqnarray*}
&&-\langle p_t(0),y_0\rangle +  \langle p(0),y_1\rangle\cr
&=& \int_0^T\langle f, y \rangle dt+\int_0^T\langle p, BB^*y_t \rangle dt+(p_{0T},y_t(T))-(p_{1T},y(T))\cr
&\leq& \int_0^T\langle f, y\rangle dt+\int_0^T\langle p, BB^*y_t \rangle dt+c(T+1)^{-\frac{1}{2}k\varrho}\|(p_{0T}, p_{1T})\|_{\mathcal{H}'} \|(y_0, y_1)\|_{\mathcal{D}(\mathcal{A}^k)}\cr
&\leq& \tilde{c} (g_1(T))^{\frac{1}{2}} [( \int_0^T \|B^*p\|_{H}^2 dt)^{\frac{1}{2}}
\!+\!(\int_0^T \|f\|_{X_{-\!{1}\!/\!{2}}}^2)^{\frac{1}{2}}+\|(p_{0T}, p_{1T})\|_{\mathcal{H}'}] \|(y_0, y_1)\|_{\mathcal{D}(\mathcal{A}^k)}.
\end{eqnarray*}

By Riesz representation theorem,  we can choose $(y_0, y_1)\in \mathcal{D}(\mathcal{A}^{k})$ such that $ \| (y_0, y_1)\|_{\mathcal{D}(\mathcal{A}^{k})}=\|(p(0), p_t(0))\|_{\mathcal{D}(\mathcal{A}^{k})'}$ and  $$-\langle p_t(0),y_0\rangle +\langle p(0),y_1\rangle=\|(p(0), p_t(0))\|_{\mathcal{D}(\mathcal{A}^{k})'}^2.$$
Thus, 
\begin{equation}
\|(p(0), p_t(0))\|_{\mathcal{D}(\mathcal{A}^{k})'}^2\leq \tilde{c} g_1(T)[ \int_0^T\|B^* p\|_{H}^2dt+\int_0^T\|f\|_{X_{-{1}\!/\!{2}}}^2dt+ \|(p_{0T}, p_{1T})\|_{\mathcal{H}'}^2]
\end{equation}
	where  $k>0$ and  $\mathcal{H}'=H\times X_{-{1}\!/\!{2}}$ is the dual space of $\mathcal{H}$. 
	The proof is completed.\hfill$\Box$

By the similar discussion, together with (H2), we obtain that	\begin{lemma}\label{L-4-2}
	Assume that (H2) holds true. Then there exists a constant $\widetilde{c}>0$ satisfying
		\begin{equation}\label{nonc}
	 \|(w(T), w_t(T))\|_{\mathcal{D}(\mathcal{A}^{{k}})'}^2\!\leq\! \tilde{c}g_2(T)[ \int_0^T\|C w\|_{H}^2dt+\int_0^T\|f\|_{X_{-{1}\!/\!{2}}}^2dt+ \|(w_{0}, w_{1})\|_{{H}\!\times\! X_{-{1}\!/\!{2}}}^2].
		\end{equation}
		where {${k}\geq 1$}, $g_2(T)= \left\{
		\begin{array}{ll} \frac{(T+1)^{-{k}\eta+1}-1}{-{k}\eta+1}& {k}\eta\neq1,\\
		\ln(T+1)& {k}\eta=1,
		\end{array}
		\right.$  and $w$ is any solution to the following inhomogeneous system
		\begin{equation}\label{109+}
		\left\{
		\begin{array}{l}
		w_{tt}+A w=f,\\
		w(0)=w_0,\quad w_t(0)=w_1.
		\end{array}
		\right.
		\end{equation}
		\end{lemma}
	{\bf Proof.~} 
Similar to the proof for Lemma \ref{L-1.1}, by duality we have
 \begin{equation}\label{413}
 \langle w_t,\phi\rangle\mathop{|}\limits_0^T -  \langle w,\phi_t\rangle\mathop{|}\limits_0^T
 +
 \int_0^T \langle w,  \phi_{tt}+A\phi \rangle \,dt=
 \int_0^T \langle f,\phi\rangle\, dt
 \end{equation}
{where $\phi$ is the solution to equation (\ref{02+}).
 Using H\"{o}lder's inequality, we get
 \begin{eqnarray*}
 	| \int_0^T \langle w,  \phi_{tt}+A\phi \rangle dt| 
 	&=&| \int_0^T \langle w,  C^*C\phi_t \rangle dt| \cr
 	&\leq& 
 	c(\int_0^T\|Cw\|^2_{H}dt)^{\frac{1}{2}}(\int_0^T \|C \phi_t\|^2_{H})^{\frac{1}{2}}.
 \end{eqnarray*}

 Then, due to Lemma \ref{le1-12} and (\ref{214}),  we get that  there exists some constant $c>0$ such that
 \begin{eqnarray}\label{407+}
 &&| \int_0^T \langle w,  \phi_{tt}+A\phi \rangle \,dt|\cr 
 &\leq& c (\int_0^T \| Cw\|_{H}^2\, dt)^{\frac12} (\frac{1}{2}\|(\phi_0^T, \phi_1^T)\|^2_{\mathcal{H}}+\frac{1}{2}c_1 (T+1)^{-{k} \eta }\|(\phi_0^T, \phi_1^T)\|_{\mathcal{D}(\mathcal{A}^{{k}})}^2)^{\frac{1}{2}}\cr
 &\leq&{c}  (\frac{1}{2}+ \frac{1}{2}c_1 (T+1)^{-{k} \eta })^{\frac{1}{2}}  (\int_0^T \| Cw\|_{H}^2\, dt)^{\frac12} \|(\phi_0^T, \phi_1^T)\|_{\mathcal{D}(\mathcal{A}^{{k}})},\quad k\geq 1.
 \end{eqnarray}
 
 We also have 
 \begin{eqnarray}\label{408+}
 \int_0^T \langle f,\phi\rangle \, dt
 &\leq&  (\int_0^T \|f\|_{X_{-{1}\!/\!{2}}}^2dt)^{\frac{1}{2}} (\int_0^T \|\phi\|_{X_{{1}\!/\!{2}}}^2dt)^{\frac{1}{2}}\cr
 &\leq& c (\int_0^T \|f\|_{X_{-{1}\!/\!{2}}}^2dt)^{\frac{1}{2}} (\int _0^T(T-t+1)^{-{{k}\eta}} \|(\phi_0^T, \phi_1^T)\|_{\mathcal{D}(\mathcal{A}^{{k}})}^2dt)^\frac{1}{2}\cr
 &\leq&{c}(g_2(T))^{\frac{1}{2}}  (\int_0^T \|f\|_{X_{-{1}\!/\!{2}}}^2dt)^{\frac{1}{2}}  \|(\phi_0^T, \phi_1^T)\|_{\mathcal{D}(\mathcal{A}^{{k}})},\; k\geq 1.
 \end{eqnarray}
 
 Hence, by (\ref{413}), (\ref{407+}) and (\ref{408+}), along with Lemma \ref{le1-12},  we obtain
 \begin{eqnarray*}
 	&&\langle w_t(T),\phi_0^T\rangle -  \langle w(T),\phi_1^T\rangle\cr
 	&=& \int_0^T\langle f, \phi \rangle dt-\int_0^T\langle w, C^*C\phi_t \rangle dt+\langle w_t(0),\phi(0)\rangle-\langle w(0), \phi_t(0)\rangle\cr
 	&\leq& \int_0^T\langle f, \phi\rangle dt-\int_0^T\langle w, C^*C\phi_t \rangle dt+c(T+1)^{-\frac{1}{2}{{k}}\eta}\|(\phi_{0}^T, \phi_{1}^T)\|_{\mathcal{D}(\mathcal{A}^{{k}})} \|(w_0, w_1)\|_{H\!\times\! X_{-{1}\!/\!{2}}}\cr
 	&\leq& {c}(g_2(T))^{\frac{1}{2}} [( \int_0^T \|Cw\|_{H}^2 dt)^{\frac{1}{2}}
 	\!+\!(\int_0^T \|f\|_{X_{-{1}\!/\!{2}}}^2)^{\frac{1}{2}}+\|(w_0, w_1)\|_{H\!\times\! X_{-{1}\!/\!{2}}}] \|(\phi_0^T, \phi_1^T)\|_{\mathcal{D}(\mathcal{A}^{{k}})}.
 \end{eqnarray*}
 
 By Riesz representation theorem,  we can choose $(\phi_0^T, \phi_1^T)\in \mathcal{D}(\mathcal{A}^{{k}})$ such that $ \| (\phi_0^T, \phi_1^T)\|_{\mathcal{D}(\mathcal{A}^{{k}})}=\|(w(T), w_t(T))\|_{\mathcal{D}(\mathcal{A}^{{k}})'}$ and  $$\langle w_t(T),\phi_0^T\rangle -  \langle w(T),\phi_1^T\rangle=\|(w(T), w_t(T))\|_{\mathcal{D}(\mathcal{A}^{{k}})'}^2.$$
 Thus, 
 \begin{eqnarray*}
&& \|(w(T), w_t(T))\|_{\mathcal{D}(\mathcal{A}^{{k}})'}^2\cr
 &\leq& \tilde{c}g_2(T)[ \int_0^T\|C w\|_{H}^2dt+\int_0^T\|f\|_{X_{-{1}\!/\!{2}}}^2dt+ \|(w_{0}, w_{1})\|_{H\!\times\! X_{-{1}\!/\!{2}}}^2],\;\;
 k\geq 1.
\end{eqnarray*}
 The proof is completed.\hfill$\Box$}

{\bf Proof of Theorem \ref{t-4-3}.~}Taking the dual product of the first and second equations in (\ref{403}) with $p^T-\bar{p}$ and $w^T-\bar{w}$, respectively, we obtain 
\begin{equation}
\int_0^T \langle(p^T-\bar{p})_{tt}, w^T-\bar{w}\rangle dt+\int_0^T\langle A (p^T-\bar{p}), w^T-\bar{w}\rangle dt=
\int_0^T\langle C^*C(w^T-\bar{w}), w^T-\bar{w}\rangle dt,
\end{equation}
and
\begin{equation}
\int_0^T\langle (w^T-\bar{w})_{tt}, p^T-\bar{p}\rangle dt +\int_0^T\langle A(w^T-\bar{w}), p-\bar{p}\rangle dt =\int_0^T\langle B(u^T-\bar{u}),p^T-\bar{p}\rangle dt.
\end{equation}
Integrating the above by parts, we obtain
\begin{eqnarray}\label{416+}
&&\int_0^T\|C(w^T-\bar{w})\|_{H}^2dt-\int_0^T\langle B(u^T-\bar{u}),p^T-\bar{p}\rangle dt\cr
&=&\langle (p^T-\bar{p})_t(T),(w^T-\bar{w})(T)\rangle-\langle(p^T-\bar{p})_t(0),(w^T-\bar{w})(0)\rangle \cr
&&-\langle (p^T-\bar{p})(T),(w^T-\bar{w})_t(T)\rangle 
+\langle (p^T-\bar{p})(0),(w^T-\bar{w})_t(0)\rangle\cr
&=&-\langle(p^T-\bar{p})_t(0),(w_0-\bar{w})\rangle +\langle \bar{p},(w^T-\bar{w})_t(T)\rangle +\langle (p^T-\bar{p})(0),w_1\rangle.
\end{eqnarray}

By Lemma \ref{L-1.1} along with (\ref{403}), we get
	\begin{eqnarray}\label{418}
&&\|(p^T(0)-\bar{p}, (p^T_t(0)-\bar{p})_t)\|_{\mathcal{D}(\mathcal{A}^{k})'}^2\cr 
&\leq&
\widetilde{c}g_1(T)[\int_0^T\|B^*(p^T-\bar{p})\|_{H}^2dt+\int_0^T\|C^*C (w^T-\bar{w})\|_{X_{-{1}\!/\!{2}}}^2dt
+\|\bar{p}\|_{H}^2].
\end{eqnarray}

Similarly, by Lemma \ref{L-4-2} along with (\ref{403}), we have 
	\begin{eqnarray}\label{419}
&&\|(w(T)-\bar{w}, w_t(T))\|_{X_{-{{k}}\!/\!{2}}\times X_{-{({k}+1)}\!/\!{2}}}^2 \cr
	&\leq&
	\widetilde{c} g_2(T)[\int_0^T\|B(u^T\!-\!\bar{u})\|_{X_{-{1}\!/\!{2}}}^2 dt\!+\!\int_0^T\|C(w^T-\bar{w})\|_{H}^2dt\cr
	&&\!+\!\|w_{0}\!-\!\bar{w}\|_{H}^2\!+\!\|w_{1}\|_{X_{-{1}\!/\!{2}}}^2],\; k\geq 1.
\end{eqnarray}

Thus, by (\ref{418}), it holds that
\begin{eqnarray}\label{421}
&&\langle(p^T-\bar{p})_t(0),(w_0-\bar{w})\rangle+\langle (p^T-\bar{p})(0),w_1\rangle\cr
&\leq& \langle \|(p^T-\bar{p})(0), (p^T-\bar{p})_t(0))\|_{\mathcal{D}(\mathcal{A}^{k})'}   \|(w_0-\bar{w}, w_1)\|_{\mathcal{D}(\mathcal{A}^{k})}\cr 
&\leq &\widetilde{c}(g_1(T))^{\frac{1}{2}}[\int_0^T\|B^*(p^T-\bar{p})\|_{H}^2dt+\int_0^T\|C^*C (w^T-\bar{w})\|_{X_{-{1}\!/\!{2}}}^2dt
+\|\bar{p}\|_{H}^2]^{\frac{1}{2}}\cr
&&\cdot \| (w_0-\bar{w}, w_1) \|_{\mathcal{D}(\mathcal{A}^k)},
\end{eqnarray}
and by  (\ref{419}), we get
\begin{eqnarray}\label{422}
&&\langle \bar{p},(w^T-\bar{w})_t(T)\rangle \leq
\| \bar{p}\|_{X_{{({k}+1)}\!/\!{2}}} \|(w^T-\bar{w})_t(T)\|_{X_{-{({k}+1)}\!/\!{2}}}\cr
&\leq& 	\widetilde{c} (g_2(T))^{\frac{1}{2}}[\int_0^T\!\|B(u^T\!-\!\bar{u})\|_{X_{-{1}\!/\!{2}}}^2 dt\!+\!\int_0^T\!\|C(w^T\!-\!\bar{w})\|_{H}^2dt
\!+\!\|w_{0}\!-\!\bar{w}\|_{H}^2\!+\!\|w_{1}\|_{X_{-{1}\!/\!{2}}}^2]^{\frac{1}{2}}\cr
&&\cdot \| \bar{p}\|_{X_{{({k}+1)}\!/\!{2}}},\; k\geq 1.
\end{eqnarray}

Substituting (\ref{421}) and (\ref{422}) into (\ref{416+}) yields

\begin{eqnarray}\label{420}
&&\int_0^T\|C(w^T-\bar{w})\|_{H}^2dt-\int_0^T\langle B(u^T-\bar{u}),p^T-\bar{p}\rangle dt\cr
&\leq&\widetilde{c}(g_1(T))^{\frac{1}{2}}[\int_0^T\|B^*(p^T-\bar{p})\|_{H}^2dt+\int_0^T\|C^*C (w^T-\bar{w})\|_{X_{-{1}\!/\!{2}}}^2dt
+\|\bar{p}\|_{H}^2]^{\frac{1}{2}}\cr
&&\cdot \| (w_0-\bar{w}, w_1) \|_{\mathcal{D}(\mathcal{A}^k)}\cr
&&+\widetilde{c} (g_2(T))^{\frac{1}{2}}[\int_0^T\!\|B(u^T\!-\!\bar{u})\|_{X_{-{1}\!/\!{2}}}^2 dt\!+\!\int_0^T\!\|C(w^T\!-\!\bar{w})\|_{H}^2dt
\!+\!\|w_{0}\!-\!\bar{w}\|_{H}^2\!+\!\|w_{1}\|_{X_{-{1}\!/\!{2}}}^2]^{\frac{1}{2}}\cr
&&\cdot \|\bar{p}\|_{X_{{(\widetilde{k}+1)}\!/\!{2}}},\quad k>0,\; \widetilde{k}\geq 1,
\end{eqnarray}
and hence, due to $C^*\in \mathcal{L}(H, X_{-{1}\!/\!{2}})$ and the third equation in (\ref{403}), we get that there exist some constants $c_j, j=1,2$ such that 
\begin{eqnarray}
&&\int_0^T\|C(w^T-\bar{w})\|_{H}^2dt-\int_0^T\langle B(u^T-\bar{u}),p^T-\bar{p}\rangle dt\cr
&\leq&c_1 g_1(T) \| (w_0-\bar{w}, w_1) \|_{\mathcal{D}(\mathcal{A}^k)}^2+c_2 g_2(T) \|\bar{p}\|_{X_{{(\widetilde{k}+1)}\!/\!{2}}}^2
\end{eqnarray}
and so,
\begin{eqnarray}\label{340}
&&\frac{1}{T}\int_0^T\|C(w^T-\bar{w})\|_{H}^2dt-\frac{1}{T}\int_0^T\langle B(u^T-\bar{u}),p^T-\bar{p}\rangle dt \cr
&\leq & \frac{c_1 g_1(T)}{T} \| (w_0\!-\!\bar{w}, w_1) \|_{\mathcal{D}(\mathcal{A}^k)}^2+\frac{c_2 g_2(T)}{T}\|\bar{p}\|_{X_{{(\widetilde{k}\!+\!1)}\!/\!{2}}}^2\!\to\! 0,\; as\; k>0, \widetilde{k}\geq 1.
\end{eqnarray}

Integrating the first equation in (\ref{403}) from $0$ to $T$ yields
$$\int_0^T A(w-\bar{w})dt=\int_0^T B(u-\bar{u})dt-w_t(T)+w_1.$$

 So,
 $$\|\frac{1}{T}\int_0^T A(w-\bar{w})dt\|_{X_{-{1}\!/\!{2}}}\leq
 \|\frac{1}{T}\int_0^T B(u-\bar{u})dt\|_{X_{-{1}\!/\!{2}}}+\frac{1}{T}(\|w_t(T)\|_{X_{-{1}\!/\!{2}}}+\|w_1\|_{X_{-{1}\!/\!{2}}}).
 $$
 
 Thus, 
 \begin{eqnarray}
 &&\|\frac{1}{T}\int_0^T (w-\bar{w})dt\|_{X_{{1}\!/\!{2}}}\cr
 &\sim& \|\frac{1}{T}\int_0^T A(w-\bar{w})dt\|_{X_{-{1}\!/\!{2}}}\cr
 &\leq&
 \|\frac{1}{T}\int_0^T B(u-\bar{u})dt\|_{X_{-{1}\!/\!{2}}}+\frac{1}{T}(\|w_t(T)\|_{X_{-{1}\!/\!{2}}}+\|w_1\|_{X_{-{1}\!/\!{2}}}).
 \end{eqnarray}
 By (\ref{340}), along with the H\"{o}lder's inequality, we get 
 \begin{eqnarray*}
  &&\|\frac{1}{T}\int_0^T (w-\bar{w})dt\|_{X_{{1}\!/\!{2}}}^2\cr
  &\leq& 
 \frac{1}{T}\int_0^T \|(u-\bar{u})\|_H^2dt+\frac{c}{T^2}\cr
  &\leq& \frac{c_1g_1(T)}{T} \| (w_0-\bar{w}, w_1) \|_{\mathcal{D}(\mathcal{A}^k)}^2+\frac{c_2g_2(T)}{T}\|\bar{p}\|_{X_{{(\widetilde{k}+1)}\!/\!{2}}}^2+ \frac{c}{T^2}\to 0,\; as \; k>0,\; \widetilde{k}\geq 1.
 \end{eqnarray*}
The proof is completed.\hfill$\Box$

\section{Examples}\label{examples}
This section is devoted to presenting some examples on hyperbolic LQ optimal control problems with weak controllability and observability. By using the abstract results obtained  in this work, we can identify the explicit slow decay rates and turnpike property for these examples.
\subsection{Wave networks}

Wave networks have been studied by many researchers (see \cite{dager}, 	 \cite{hanzuazua}, \cite{lagnese94}, \cite{2}, \cite{xusiam08} and the references therein). Here we consider a simple case: star-shaped wave networks.

 The setting as given in \cite{dager} is chosen to form the wave networks.
A 1-d wave 
network $\mathcal{R}:= \cup_{j=1}^N e_j$  is formed by $N$ wave equations on the curve $e_j,\; j=1,2\cdots,N$ with interval $(0,\ell_j)$.  For $k\neq j$, $\overline{e_j}\cap \overline{e_k}$, where $\overline{e_i},\; i=1,2,\cdots, N$ is denoted by the closure of $e_i$, is either empty or a common end called a vertex or a node. Assume that  the wave
equation arises on the intervals $(0, \ell_j),\; j=1, 2,\cdots,
N$ in the network with state $(w_j, w_{j,t})$, respectively.

Let $G=\{e_j,\; j=1,2,\cdots N\}$ be the set of $e_j$ of $\mathcal{R}$,  and  $\mathcal{V}$ be
the set of vertices of $\mathcal{R}$. Denote by $G_v=\{ j=\{1,2,\cdots, N\},\; v\in \overline{e_j}\}$
the set of edges having $v$ as a vertex. Denote by $card(G_v)$ the number of edges that meet at $v$. We call $v$ is an exterior node if $card(G_v)=1$, the set of which is denoted by $\mathcal{V}_{int}$, while if $card (G_v)\geq 2$, the node $v$ is called an interior node and the set of them is $\mathcal{V}_{ext}$.

Assume that the Dirichlet conditions are fulfilled at the exterior
nodes and the geometrical continuity is satisfied at the interior nodes
of the network. The control is assumed to be located only at one edge. For convenience, we set the index of the controlled edge is $1$ and $j_0$ is the index of the observed edge.   Then we get  the following wave equations on
a network (see Fig. \ref{fff2} for instance):
\begin{figure}
	\begin{center}
		\includegraphics[width=8cm]{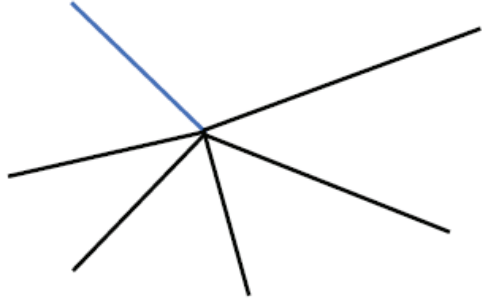}
%		\put(-70,-9){$\ell_1$} \put(-71,30){$e_1$} \put(-40,40){$\ell_2$}
%		\put(-65,65){$e_2$} \put(-40,70){$\vdots$} \put(0,110){$\ell_{N_1}$}
%		\put(-50,99){$e_{N_1}$} \put(-100,140){$\ell_{N_1\!+\!1}$}
%		\put(-115,105){$e_{N_1\!+\!1}$} \put(-180,115){$\ell_{N_1\!+\!2}$}
%		\put(-158,85){$e_{N_1\!+\!2}$} \put(-158,70){$\vdots$}
%		\put(-183,60){$\ell_{N\!-\!1}$} \put(-136,55){$e_{N\!-\!1}$}
%		\put(-140,5){$\ell_N$}
%		\put(-113,30){$e_N$}
%		\put(-83,77){$0$}
		\caption{Star-shaped network of wave equations}\label{fff2}
	\end{center}
\end{figure}
\begin{equation}\label{20001+}
\left\{\begin{array}{l}
w_{1,tt}(x,t)- w_{1,xx}(x,t)=u(x,t),\; x\in(0, \ell_1),\; t>0,\\
w_{j,tt}(x,t)- w_{j,xx}(x,t)=0,\; x\in(0, \ell_j),\; j=2,...N,\; t>0,\\
w_\ell(v,t)=w_j(v,t),\; \forall \ell, j\in G_v, v\in \mathcal{V}_{int},\; t>0,\\
\sum_{j\in G_v}\limits
w_{j,x}(v,t)=0,\; \forall v\in \mathcal{V}_{int},\; t>0,\\
w_j(v,t)=0,\; \forall j\in G_v, v\in \mathcal{V}_{ext},\; t>0,\\
w_j(t=0)=w_j^{0},\; w_{j,t}(t=0)=w_j^{1}, \quad j=1,
2,\cdots, N,
\end{array}
\right.
\end{equation}
where $( (w_j^{0})_{j={1}}^N,
(w_j^{1})_{j={1}}^N)$ is the given initial state.

Consider the  infinite-horizon optimal control problem of quadratic type
\begin{equation}\label{optinf-network}
\min J^\infty(u)=\frac{1}{2}\int_0^\infty [\|u(x,t)\|^2_{U}+ \|w_{j_0,x}(x,t)\|_{L^2(0,\ell_{j_0})}^2]dt,\; u\in L^2(0,T; U),
\end{equation}

We define the Hilbert space
$$L^2(\mathcal{R}) =
		\{f=(f_j)_{j=1}^{N}| f_j\in
		L^2(0,\ell_j),\; \forall j=1, 2,\cdots, N\}$$
		and
		$$
		V^m(\mathcal{R})=\{f=(f_j)_{j=1}^{N}|f_j\in H^m(0,\ell_j),
		f_j(\ell_j)=0, f_j(0)=f_i(0), \forall i, j=1,2,\cdots,N\}.$$
		
Define the operator $A$ in $L^2(\mathcal{R})$ as 
$$A(w_j)_{j=1}^N := -(w_{j,xx})_{j=1}^N$$
with domain
$\mathcal{D}({A})\!=\!\left\{ (w_j)_{j=1}^N\in {V}^2(\mathcal{R}) \left|
\sum_{j=1}^N\limits w_{j,x}(0)=0
\right.\right\}.$ Thus,  the star-shaped networks (\ref{20001+}) can be rewritten as the abstract form (\ref{001}).

Set the state space  $\mathcal{H}$ as follows:
$$\mathcal{H}=V^1(\mathcal{R})\times L^2(\mathcal{R})$$ equipped with inner product: for
$W=((w_j)_{j=1}^N, (z_j)_{j=1}^N),\; \widetilde{W}=((\widetilde{w})_{j=1}^N, (\widetilde{z})_{j=1}^N) \in \mathcal{H}$,
$$(W,\widetilde{W})_\mathcal{H}= \sum_{k=1}^{N}\int_0^{\ell_j}w_{k,x}\overline{\widetilde{w}_{j,x}}dx
+
\sum_{j=1}^N \int_0^{\ell_j} z_{j}\overline{\widetilde{z}_{j}}dx.$$

It is easy to check that $(\mathcal{H}, \|\cdot\|_{\mathcal{H}})$ is
a Hilbert space. The system operator can be set as $\mathcal{A}:=\left[ \begin{array}{cc} 0&I \cr -A& 0\end{array} \right]$ in $\mathcal{H}$. 

By \cite{dager}, we have the following weak observability estimate for system ({\ref{20001+}}).
\begin{proposition}\label{hp-5-1}
	There exists a positive constant ${T}>0$  such that
	\begin{equation}\label{0405+++}
	c_1 \int_0^T\int_0^{\ell_1} |w_{1,t}|^2dxdt\geq \sum_{n\geq 1}\limits \gamma_n^2[\lambda_n^2a_n^2+b_n^2]
	\end{equation}
	and 
	\begin{equation}
		c_2 \int_0^T\int_0^{\ell_{j_0}} |w_{j_0,x}|^2dxdt\geq \sum_{n\geq 1}\limits \widetilde{\gamma}_n^2[\lambda_n^2a_n^2+b_n^2] 
		\end{equation}
	where  $c_j>0,\; j=1,2 $ are some constant, $\lambda_n^2$ is the eigenvalue of the operator $A$ corresponding to system (\ref{20001+}),  $a_n,\; b_n$ are the Fourier coefficients given as in (\ref{fourier}), and $\gamma_n^2,\; \widetilde{\gamma}_n^2>0$ are the weights,  which are determined by the lengths of edges involved in the network.  
\end{proposition}

In general, we just know that $\gamma_n\to 0$ as $n$ goes to infinity and can not get a better estimate for the decay rate of $\gamma_n$. However, following the proof for Theorem \ref{th1}, as well as Remark \ref{remark1}, if the estimates $\gamma_n\geq \Phi_1(\lambda_n)$ and $\widetilde{\gamma}_n\geq \Phi_2(\lambda_n)$ hold,  we can get 
	\begin{eqnarray}\label{lowupper-network+}
&&\widetilde{c}_1\|( (w_j^{0})_{j={1}}^N,
(w_j^{1})_{j={1}}^N)\|_{{\mathcal{D}(\!\Phi_2(\!\mathcal{A}\!)\!)}}^2\leq  \langle\hat{E} ( (w_j^{0})_{j={1}}^N,
(w_j^{1})_{j={1}}^N),( (w_j^{0})_{j={1}}^N,
(w_j^{1})_{j={1}}^N)\rangle,\\
\cr
&&\langle\hat{E}(\!(w_j^{0})_{j={1}}^N,
(w_j^{1})_{j={1}}^N\!),(\!(w_j^{0})_{j={1}}^N,
(w_j^{1})_{j={1}}^N\!)\rangle\leq c_2\|( (w_j^{0})_{j={1}}^N,
(w_j^{1})_{j={1}}^N)\|^2_{\mathcal{D}(\!\Phi_1^{-1}(\mathcal{A})\!)\!},
\end{eqnarray}
where  $c_1, c_2$ are some positive constants, $\mathcal{D}(\!\Phi_1^{-1}(\mathcal{A})\!)$ is denoted by the space satisfying 
$$ \sum_{n\geq 1}\limits \Phi_1^{-2}(\lambda_n)[\lambda_n^2a_n^2+b_n^2] < \infty,$$
and  $a_n,\; b_n$ are the Fourier coefficients given as in (\ref{fourier}).
%Thus, according to the proof in Theorem \ref{th1}, we get the following slow decay rate of the closed-loop system with optimal feedback control.
%	$$ \| ((\hat{w_j})_{j=1}^N, (\hat{w}_{j,t})_{j=1}^N)\|_{\mathcal{H}}^2\leq {M}{(t+1)^{-1}} \|( (w_j^{0})_{j={1}}^N,
%(w_j^{1})_{j={1}}^N)\|_{\mathcal{D}(\! \frac{1}{\Phi(\mathcal{A})})\!}^2.$$

	Although it is unknown on the estimate of $\gamma_n$ for general networks, it can be better estimated for some special cases. For instance,
	by \cite{dager}, we  see that  for star-shaped networks, the weights $\gamma_n$  is determined by
	the edge-length ratios $\frac{\ell_i}{\ell_j}$, where $i,j=
	2,\cdots,N,\; i\neq j$ and $\widetilde{\gamma}_n$ is determined by the ratios $\frac{\ell_i}{\ell_j}$, where $i,j=
	1, 2,\cdots,N,\; i\neq j,\; i\neq j_0,\; j\neq j_0$.
%	 and $\widetilde{\gamma}_n$ in
%	(\ref{0405+++}) are given as follows
%	$$\gamma_n=\max_{i=2,\cdots,N}\prod_{ j\neq i}\limits|\sin(\lambda_n\ell_j)|,\; \forall n\geq 1,$$
% and   condition $\inf_{n>0} \limits\gamma_n^2=c>0$, which is corresponding to the exact observability, does not  hold
%	for the networks. Indeed,  it  always satisfies   $\liminf_{n\to\infty}
%	\limits\gamma_n^2=0$ and  $\gamma_n^2$ is determined by
%	the ratios $\frac{\ell_i}{\ell_j}$, where $\ell_i, \ell_j,\; i,j=
%2,\cdots,N,\; i\neq j$. 
Specifically, if $\frac{\ell_i}{\ell_j}$
 belongs to some special irrational sets(see \cite{Schmidt}), we can obtain that  there always exist some constants $\zeta, \xi>0$ and $c>0$ such that $\gamma_n\geq \frac{c}{\lambda_n^\zeta}$ and $\gamma_n\geq \frac{c}{\lambda_n^{\xi}}$. Based on it, together with Theorem \ref{th1},
 we can obtain the following slow decay rate of the closed-loop system with  Riccati-based optimal feedback control. 
 $$ \| ((\hat{w_j})_{j=1}^N, (\hat{w}_{jt})_{j=1}^N)\|_{\mathcal{H}}^2\leq {M}{(t+1)^{-\frac{s}{\zeta+\xi}}} \|( (w_j^{0})_{j={1}}^N,
 (w_j^{1})_{j={1}}^N)\|_{\mathcal{D}(\mathcal{A}^{\zeta+\xi+s})}^2, \quad  s>0. $$
Besides,   if $((w_j^0-\bar{w}_j)_{j=1}^N, (w_j^1)_{j=1}^N)\in \mathcal{D}(\mathcal{A}^k),\; k>0$ and $(\bar{p}_j)_{j=1}^N\in V_{2}(\mathcal{R})$ hold, the averaged turnpike property (\ref{281}), (\ref{282}) hold for such kind of networks,
that is,
$$
\|\frac{1}{T}\int_0^T(u^T-\bar{u})dt\|_{L^2(0,\; \ell_1)}\to 0, \; as\;\;  T\to \infty
$$
and
$$
\|\frac{1}{T}\int_0^T ((w_j^T(t))_{j=1}^N-(\bar{w}_j)_{j=1}^N)dt\|_{V^1(\mathcal{R})}\to 0,  \; as\;\;  T\to \infty.
$$

\subsection{Wave equation without GCC: Rectangular domain}

Consider a wave equation on a square $\Omega=(0,\pi)\times (0,\pi)$ with local control(see Fig. \ref{tank}):
\begin{figure}
	\begin{center}
		\includegraphics[width=8cm]{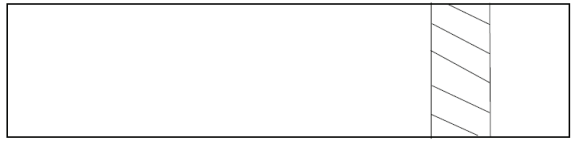}
		%		\put(-70,-9){$\ell_1$} \put(-71,30){$e_1$} \put(-40,40){$\ell_2$}
		%		\put(-65,65){$e_2$} \put(-40,70){$\vdots$} \put(0,110){$\ell_{N_1}$}
		%		\put(-50,99){$e_{N_1}$} \put(-100,140){$\ell_{N_1\!+\!1}$}
		%		\put(-115,105){$e_{N_1\!+\!1}$} \put(-180,115){$\ell_{N_1\!+\!2}$}
		%		\put(-158,85){$e_{N_1\!+\!2}$} \put(-158,70){$\vdots$}
		%		\put(-183,60){$\ell_{N\!-\!1}$} \put(-136,55){$e_{N\!-\!1}$}
		%		\put(-140,5){$\ell_N$}
		%		\put(-113,30){$e_N$}
		%		\put(-83,77){$0$}
		\caption{ Wave equation on rectangular domain with local control}\label{tank}
	\end{center}
\end{figure}
\begin{equation}\label{710}
\left\{
\begin{array}{l}
w_{tt}=\Delta w(x,t)- \chi_{{\Omega}_0} u(x,t),\\
w|_{\partial \Omega}=0,\\
w(x,0)=w_0(x),\; w_t(x,0)=w_1,
\end{array}
\right.
\end{equation}
where $\Omega_0=\{(x_1,x_2)| a<x_1<b,\; 0<x_2<\pi\}$ is a strip-type subdomain parallel to one boundary. Meanwhile,  assume that  $$Cw=\nabla w  \; in\;  \Omega,$$ Obviously, $(A,C)$ is exactly observable.

The stability of wave equation in rectangular domain with locally viscous damping  was ever considered by \cite{batty} and \cite{rein}, where  they obtained that  under $u(x,t)=\alpha w_t(x,t)$,
 the system can be stabilized polynomially,  and the optimal decay rate is given as follows:
 \begin{equation} \label{706}
 \|(w,w_t)\|^2_\mathcal{H}\leq C_k (t+1)^{-\frac{4}{3}k}\|(w_0,w_1)\|_{\mathcal{D}(\mathcal{A}^k)}^2,\quad k>0,
 \end{equation}
 in which the state space is chosen as  $\mathcal{H}=H^1_0(\Omega)\times L^2(\Omega)$ and $\mathcal{A}=\left[\begin{array}{cc}
 0&I\\
 -\Delta &0
 \end{array}\right]$.

Choose the following infinite-horizon quadratic cost performance index:
\begin{equation}\label{optinf-ran}
\min J^\infty(u)=\frac{1}{2}\int_0^\infty [\|u(x, t)\|^2_{L^2(\Omega_0)}+\|\nabla w(x,t)\|_{L^2(\Omega)}^2]dt,\; u\in L^2(0,T; U).
\end{equation}

By (\ref{706}), together with the proof of Theorem \ref{th1}, we can get the bounds related to the corresponding algebraic Riccati operator, that is,
	there exist  constant $c_j>0,\;j=1,2$ satisfying
\begin{equation}\label{707}
c_1\|(w_0,w_1)\|_{\mathcal{H}}^2\leq  \langle\hat{E} (w_0,w_1),(w_0,w_1)\rangle \leq c_2\|(w_0, w_1)\|^2_{\mathcal{D}(\mathcal{A}^{\frac{3}{4}})}.
\end{equation}

Thus, by Theorem \ref{th1}, 
there exists a  positive constant $M>0$ such that
\begin{equation}\label{708}
 \| (\hat{w}(t), \hat{w}_t(t))\|_{\mathcal{H}}^2\leq {M}{(t+1)^{-\frac{4s}{3}}} \|(w_0,w_1)\|_{\mathcal{D}(\mathcal{A}^{\frac{3}{4}+s})}^2,\quad s>0,
 \end{equation}
where $(\hat{w}(t), \hat{w}_t(t))$ is the solution to system (\ref{710})  under the Riccati-based optimal feedback control law $u(t)=-\chi_{\Omega_0}\hat{E}(w, w_t)$.

In terms of turnpike property, by Theorem \ref{t-4-3}, we further obtain that for any $ (w_0-\bar{w}, w_1) \in\mathcal{D}(\mathcal{A}^k),\; k>0$ and  $\bar{p}\in {H^2(\Omega)}$,
the averaged turpike property holds, that is,  as $T\to\infty$,
$$
\|\frac{1}{T}\int_0^T (u^T-\bar{u})dt\|_{L^2(\Omega_0)}^2 \to 0,
\quad
\|\frac{1}{T}\int_0^T (w^T(t)-\bar{w})dt\|_{H^1(\Omega)}^2 \to 0.
$$
\begin{remark}
	In this example, although there is no weak observability estimate (H1)  being fulfilled directly, by the proof in Theorem \ref{th1}, the slow decay rate of the system with locally viscous damping given as in (\ref{706}) is enough to help us obtain the upper bound in (\ref{707}). In fact,  the slow decay rate  is ``almost" equivalent to the weak observability estimate (H1).
	
It should be noted that some more general energy decay rates for multi-dimensional wave equation on partially rectangular or torus were obtained in \cite{anapde}, \cite{burq1} and \cite{phung1}, based on which, the slow decay rates and turnpike properties for such infinite-horizon LQ optimal control problems can be also  estimated similarly from Theorem \ref{th1} and \ref{t-4-3}, respectively.
	\end{remark}

\subsection{ Wave equation without GCC: General case}
Consider the wave equation on domain $\Omega$ with local control (see Fig. \ref{abstra}).
\begin{figure}
	\begin{center}
		\includegraphics[width=8cm]{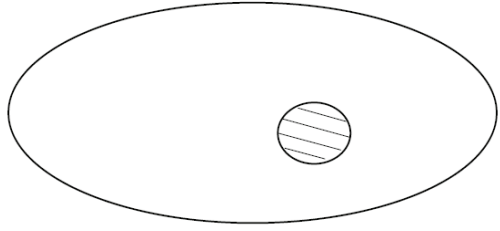}
		%		\put(-70,-9){$\ell_1$} \put(-71,30){$e_1$} \put(-40,40){$\ell_2$}
		%		\put(-65,65){$e_2$} \put(-40,70){$\vdots$} \put(0,110){$\ell_{N_1}$}
		%		\put(-50,99){$e_{N_1}$} \put(-100,140){$\ell_{N_1\!+\!1}$}
		%		\put(-115,105){$e_{N_1\!+\!1}$} \put(-180,115){$\ell_{N_1\!+\!2}$}
		%		\put(-158,85){$e_{N_1\!+\!2}$} \put(-158,70){$\vdots$}
		%		\put(-183,60){$\ell_{N\!-\!1}$} \put(-136,55){$e_{N\!-\!1}$}
		%		\put(-140,5){$\ell_N$}
		%		\put(-113,30){$e_N$}
		%		\put(-83,77){$0$}
		\caption{ Wave equation on general domain with local control}\label{abstra}
	\end{center}
\end{figure}

\be\label{gene}
\begin{cases}
	w_{tt}- \Delta w +u\chi_{\mathcal O}=0,  & \hbox{in $(0,T)\times \Omega$},\\
	w=0, & \hbox{on $(0,T)\times\partial \Omega$},
	\\ w(0)=w_0\,,\,w_t(0)=w_1,
	& 
\end{cases}
\ee
where $u$ is the control input, $\mathcal{O}$ is a subset of the whole domain $\Omega$. Suppose that $Cw=\chi_{\mathcal{\widetilde{{O}}}}\nabla w$, where  $\widetilde{\mathcal{O}}$ is another subset of $\Omega$.

Choose the following infinite-horizon quadratic cost performance index:
\begin{equation}\label{optinf-ran+}
\min J^\infty(u)=\frac{1}{2}\int_0^\infty [\|u(x, t)\|^2_{L^2(\mathcal{O})}+\|\nabla w(x,t)\|_{L^2(\widetilde{\mathcal{O}})}^2]dt,\; u\in L^2(0,T; U). 
\end{equation}
In \cite{zuazua2},  Porretta and Zuazua showed that 
the solution to  system (\ref{gene}) with the  Riccati-based  optimal feedback control can be stabilized exponentially in the energy space $H^1(\Omega)\times L^2(\Omega)$, provided that the subset $\mathcal{O} \subset \Omega$  and $\widetilde{\mathcal{O}} \subset \Omega$ verifies the GCC. It is well-known that the GCC guarantees the exact observability of $(A, B^*)$ and $(A,C)$.

Whenever ${\mathcal{O}}$ and $\widetilde{\mathcal{O}}$ are  general open nonempty subsets of $\Omega$, not necessarily satisfying the GCC,  the exact observability of $(A, B^*)$ and $(A,C)$ cannot be fulfilled. 
However,  some weak observability estimates still hold for $(A, B^*)$ and $(A,C)$. In fact, from \cite{ZZ},  we have that  for system (\ref{gene}),  there exist some constants $a>0$, $c>0$ and $T>0$ satisfying 
$$
c\int_0^T
\int_{\mathcal{O}}|w_t|^2dxdt\geq \| (w_0, w_1)\|_{\mathcal{D}(e^{-a\mathcal{A}})}^2  \quad ((A,B^*)\; weakly\; obserbvable),$$
 and the following logarithmic decay rate holds under feedback control $u(t)=\chi_{\mathcal O}w_t$.
\begin{equation}\label{decay}
\into [|w_t(t)|^2+|\nabla w(t)|^2]dx \leq \frac{C_0}{[\log(2+t)]^2}\left( \|w_0\|^2_{H^2(\Omega)}+ \|w_1\|^2_{H^1(\Omega)}\right) \qquad \forall t>0\,.
\end{equation}

Similarly, there exist some constants $b>0$, $c>0$ and $T>0$ satisfying 
$$
c\int_0^T\int_{\widetilde{\mathcal{O}}}|\nabla w|^2dxdt\geq \| (w_0, w_1)\|_{\mathcal{D}(e^{-b\mathcal{A}})}^2, \quad ((A,C)\; weakly\; obserbvable).$$

Thus,  as it was presented in Remark \ref{remark1},   the lower and upper bounds can be derived by (\ref{exp1}). Moreover, we  can get the slow decay rate for system (\ref{gene}) with Riccati-based optimal feedback control, as given in (\ref{exp2}).

The averaged turnpike property still holds for  any $ (w_0-\bar{w}, w_1) \in\mathcal{D}(\mathcal{A})$ and  $\bar{p}\in {H^2(\Omega)}$.
Indeed, based on (\ref{decay}),  along with the proof for Lemma  \ref{L-1.1} and \ref{L-4-2},  we can see that the results in Lemma \ref{L-1.1} and \ref{L-4-2} still hold with some  $g_1(T)$ (resp. $g_2(T)$), which is determined by $\int _0^T \frac{1}{[\log(2+t)]^2}dt$ (resp. $\int _0^T \frac{1}{[\log(2+T-t)]^2}dt$) (see (\ref{408}), (\ref{408+})).

Note that  by mean value theorem of integrals, we have
\begin{equation}
\int _0^T \frac{1}{[\log(2+t)]^2}dt=\int _0^T \frac{1}{[\log(2+T-t)]^2}dt= \frac{T}{[\log(2+\alpha_T T)]^2},
\end{equation}
where $\alpha_T T\to \infty$ as $T\to\infty$. Thus, choosing 
$g_1(t)=g_2(t)=\frac{T}{[\log(2+\alpha_T T)]^2}$ and
 following the proof of Theorem \ref{t-4-3}, we  finally obtain from (\ref{340}) with $k=\widetilde{k}=1$ that 
\begin{eqnarray}\label{340+}
&&\frac{1}{T}\int_0^T\|\nabla(w^T-\bar{w})\|_{L^2(\Omega)}^2dt+\frac{1}{T}\int_0^T \|u^T-\bar{u}\|_{L^2(\mathcal{O})}^2 dt \cr
&\leq & \frac{c_1 g_1(T)}{T} \| (w_0-\bar{w}, w_1) \|_{\mathcal{D}(\mathcal{A})}^2+\frac{c_2 g_2(T)}{T}\|\bar{p}\|_{H^2(\Omega)}^2\cr
&=& \frac{1}{[\log(2+\alpha_T T)]^2}(c_1\| (w_0-\bar{w}, w_1) \|_{\mathcal{D}(\mathcal{A})}^2+{c_2}\|\bar{p}\|_{H^2(\Omega)}^2)\to 0.
\end{eqnarray}

It should be noted that the  estimate (\ref{340+}) is consistent with the one in Section 4.3 in Porretta and Zuazua \cite{zuazua2}, where the LQ optimal control problem related to the  multi-dimensional  wave equation without GCC was considered.

\section{Conclusions}\label{open}

In this work,  we considered the slow decay rate and  turnpike property  of the hyperbolic LQ optimal control problems and mainly  obtained the following results:

1. The slow decay rate of infinite-horizon hyperbolic LQ optimal control problems was considered. Under the weak observability of $(A,B^*)$ and $(A,C)$,  the lower and upper bounds of the corresponding algebraic Riccati operator were estimated, respectively. Then the explicit slow decay rate of the closed-loop system with Riccati-based optimal feedback control  was estimated, which is a key step to further discuss the turnpike properties of the LQ optimal control problems if lacking of exact controllability or observability.

 2. Under weak observability of $(A,B^*)$ and $(A,C)$ hypothesises, the averaged turnpike property for the LQ optimal control problems  was  proved under certain conditions on the regularity of the  initial state, the  stationary optimal state and its dual. This result is consistent with the  slow turnpike in average identified in section 4.3 in  \cite{zuazua2} which can be considered as one special case of our work (see  the example in Section 4.3 in our work).

Besides the averaged turnpike property,  we would like to point out that the starting point of our work was to see  whether  there were some kinds of point-wise slow turnpike properties holding for the LQ optimal control problems under weak controllability or observability hypothesises.  In fact, it can be seen from \cite{zuazua2} and \cite{zuazua3} that  the exponential decay of the closed-loop systems with Riccati-based optimal feedback control  always leads to the exponential-type point-wise turnpike property when exact controllability and observability hypothesises are fulfilled, that is, there exist $\lambda>0$ and $K>0$ such that
\be\label{wavex}
\| w^T(t)-\bar w\|_{X_{1/2}}^2+ \|u^T(t)-\bar u\|_{H}^2 \leq K (e^{-\la t}+ e^{-\la (T-t)})\qquad \forall t\in [0,T]\,,
\ee
where $(w^T, u^T)$ and $(\bar w, \bar u)$ are the optimal trajectories and controls for the optimal control problem (\ref{opt-z}) and for the stationary problem (\ref{staopt}), respectively.
Thus, based on the slow decay rate obtained in our work,  it is reasonable to predict that there should be some kinds of slow point-wise turnpike properties  holding under weak controllability or observability.

However, note that the  energy estimate method based on the properties of Riccati operator proposed in \cite{zuazua2} can not work well for this issue under consideration. 
Indeed,  when doing so, there is an inevitable process to estimate  $\|w^T-\bar{w}\|_{X_{1/2}}$ by using Duhamel's formula along with Gronwall's lemma. This process is tough to tackle  for the LQ optimal control problems under weak observability of $(A,B^*)$ and $(A,C)$  like hypothesises (H1) and (H2),  because that  under these weak hypothesises, the corresponding closed-loop systems with Riccati-based optimal controls always achieve  non-uniform slow decay rates which depending on the regularity of the initial states, 
while it can be carried out effectively for  the case of exact controllability and observability due to the uniform exponential decay rates always hold for all initial states and this  uniformity causes the  Gronwall's lemma easy to be used.
So, it is still an open problem  that whether  there are some kinds of  slow point-wise turnpike property holding for such hyperbolic LQ problems. Compared to (\ref{wavex}), 
based on the weak observation of $(A,B^*)$ and $(A,C)$ as given in (H1) and (H2),  it seems that  the slow point-wise turnpike property could hold and  have the following  form:
 $$\|w^T-\bar{w}\|_{X_{{1}\!/\!{2}}}^2+\|u^T(t)-\bar{u}\|_H^2\leq K (t+1)^{-k_1 } \|(w_0, w_1)\|_{\mathcal{D}(\mathcal{A})}^2+(T-t+1)^{-k_2}\|\bar{p}\|_{X_1}^2.$$ 
 It could be verified from the view of frequency domain by Riesz basis representation, but  careful estimates are still needed,  which is an interesting issue and worth investigating in future.
 
 The same problems are also worth discussing for  linear parabolic systems or some non-linear systems such as semilinear wave equations,  nonlinear models of fluid mechanics and so on. 
 In addition, the case of  weak boundary observability  or controllability is another interesting issue. Some new techniques could be involved due to the  unboundedness of the boundary control or observation operators.

%\section*{Acknowledgments} We warmly thank Emmanuel Tr\'elat for useful discussions on the topic of this article.  We also thank the editors and referees for fruitful suggestions, in particular concerning the  connection of our work with the turnpike property in the mathematical economics theory.
%
%A. Porretta wishes to thank the BCAM Center in Bilbao for the invitation and for providing a stimulating environment for  collaboration. 

\end{document}